\documentclass[10pt]{article}
\usepackage{amsmath, amsthm}
\usepackage{amscd}
\usepackage{amsfonts}
\usepackage{amssymb}
\usepackage{mathrsfs}
\usepackage{stmaryrd}
\usepackage{color}

\newcounter{ccomm}
\definecolor{blue}{rgb}{0,0,.6}

\makeindex

\newtheorem{theorem}{Theorem}[section]
\newtheorem{proposition}[theorem]{Proposition}
\newtheorem{lemma}[theorem]{Lemma}
\newtheorem{corollary}[theorem]{Corollary}

\theoremstyle{definition}
\newtheorem{definition}[theorem]{Definition}

\newtheorem{remark}[theorem]{Remark}

\def\proof{{\noindent\sc Proof. \quad}}
\newcommand\proofof[1]{{\noindent\sc Proof of {#1}. \quad}}
\def\eproof{{\mbox{}\hfill\qed}\medskip}
\begin{document}

\makeatletter

\renewcommand{\bar}{\overline}

%Other math symbols
\newcommand{\x}{\times}
\newcommand{\<}{\langle}
\renewcommand{\>}{\rangle}
\newcommand{\into}{\hookrightarrow}

%Greek letters
\renewcommand{\a}{\alpha}
\renewcommand{\b}{\beta}
\renewcommand{\d}{\delta}
\newcommand{\D}{\Delta}
\newcommand{\e}{\varepsilon}
\newcommand{\g}{\gamma}
\newcommand{\G}{\Gamma}
\renewcommand{\l}{\lambda}
\renewcommand{\L}{\Lambda}
\newcommand{\n}{\nabla}
\newcommand{\var}{\varphi}
\newcommand{\s}{\sigma}
\newcommand{\Sig}{\Sigma}
\renewcommand{\t}{\theta}
\renewcommand{\O}{\Omega}
\renewcommand{\o}{\omega}
\newcommand{\z}{\zeta}
\newcommand{\balpha}{\boldsymbol \alpha}
\newcommand{\ab}{\alpha_{\bullet}}

%Other macros
\newcommand{\p}{\partial}
\renewcommand{\hat}{\widehat}
\renewcommand{\bar}{\overline}
\renewcommand{\tilde}{\widetilde}
\def\Dn{\mathcal{D}}
\def\oU{\bar{U}}
\def\MD{{\sf{MD}}}
\def\MDF{{\sf{MDVar}}}
\def\MDFix{{\sf{MDFix}}}
\def\KFix{K_{\mathsf{Fix}}}
\def\bz{{\boldsymbol z}}
\def\bu{\overline{u}}

%%%%%%%%
% Some fonts
%%%%%%%

%\newcommand{\fiverm}{\tiny\rm}
\font\eightrm=cmr8
\font\ninerm=cmr9

%%%%%
% The next macros define fonts for reals, rationals, complex,
% integers and natural numbers by \R,\Q,\C,\Z and \N respectively.
% Also, \Ri gives \R to the infinity.
%%%%%

\def\N{\mathbb{N}}
\def\Z{\mathbb{Z}}
\def\R{\mathbb{R}}
\def\E{\mathop{\mathbb E}}
\def\Q{\mathbb{Q}}
\def\C{\mathbb{C}}
\def\F{\mathbb{F}}
\def\proj{\mathbb{P}}
\def\vol{\mathsf{vol}}

%%%%%%%
% For writing algorithms
%%%%%%%

\def\algo{\begin{center}
               \begin{minipage}{6in}
               \begin{tabbing}
               \marks}
\def\falgo{\end{tabbing}
                \end{minipage}
                \end{center}}
\def\marks{nn\= nn\= nn\= nn\= nn\= nn\= nn\= \kill}

%%%%%%%
% Some mathematical operators and constants
%%%%%%%

\def\ll{{[\kern-1.6pt [}}
\def\rr{{]\kern-1.4pt ]}}
\def\bll{{\biggl[\kern-3pt \biggl[}}
\def\brr{{\biggr]\kern-3pt \biggr]}}
\def\sg{{\rm sign}\,}
\def\sig{\overline{\rm sign}\,}
\def\absmax{\mathop{\underline{\rm max}}\limits}
\def\argmax{\mathop{{\rm argmax}}\limits}
\def\rank{{\rm rank}\,}
\def\Im{{\rm Im}}
\def\Real{{\rm Re}}
\def\dist{{\rm dist}}
\def\degree{{\rm degree}\,}
\def\grad{{\rm grad}\,}
\def\size{{\rm size}}
\def\supp{{\rm supp}}
\def\Id{{\rm Id}}
\def\fl{\mathop{\tt fl}}
\def\op{\mathop{\tt op}}
\def\sizem{{\rm size}_\mu}
\def\sizew{{\rm size}_{\rm h}}
\def\cost{{\mathsf{cost}}}
\def\ecost{{\mathsf{ept\_cost}}}
\def\vol{{\mathsf{vol}\;}}
\def\continue{{\mathsf{continue}}}
\def\halted{{\mathsf{halted}}}
\def\term{{\mathsf{term}}}
\def\card{{\rm card}}
\def\mod{{\rm mod}\;}
\def\trace{{\rm trace}}
\def\Prob{\mathop{\rm Prob}}
\def\length{{\rm length}}
\def\Diag{{\rm Diag}}
\def\diam{{\rm diam}}
\def\arg{{\rm arg}\,}
\def\ot{\leftarrow}
\def\transp{^{\rm T}}
\def\bL{{\bf L}}
\def\bR{{\bf R}}
\def\bfE{{\bf E}}
\def\bfB{{\mathbf B}}
\def\bfe{{\mathbf e}}
\def\bff{{\mathbf f}}
\def\bfg{{\mathbf g}}
\def\bd{{\bf d}}
\def\adj{{\rm adj}\:}
\def\cruz{\raise0.7pt\hbox{$\scriptstyle\times$}}
\def\todif{\stackrel{\scriptscriptstyle\not =}{\to}}
\def\nedif{\stackrel{\scriptscriptstyle\not=}{\nearrow}}
\def\sedif{\stackrel{\scriptscriptstyle\not=}{\searrow}}
\def\nodown{\downarrow\mkern-15.3mu{\raise1.3pt\hbox{$\scriptstyle\times$}}}
\def\noto{\to\mkern-20mu\cruz}
\def\notox{\to\mkern-23mu\cruz}
\def\alg{{\hbox{\scriptsize\rm alg}}}
\def\bS{\mathbb{S}}
\def\dpr{{d_{\proj}}}
\def\dsp{{d_{\bS}}}
\def\ttl{{\tt l}}
\def\tto{{\tt o}}
\def\ttml{\bar{\tt l}}
\def\ttmo{\bar{\tt o}}
\def\ri{R^\infty}
\def\Oh{{\cal O}}
\def\parts{{\cal P}}
\def\coeff{{\hbox{\rm coeff}}}
\def\Error{{\hbox{\tt Error}}\,}
\def\Bound{{\hbox{\tt Bound}}\,}
\def\RelErr{{\hbox{\tt RelError}}\,}
\def\appleq{\hbox{\lower3.5pt\hbox{$\;\:\stackrel{\textstyle<}{\sim}\;\:$}}}
\def\ALH{{\sf{ALH}}}
\def\ALHFix{{\sf{ALHFix}}}
\def\ALHF{{\sf{ALHVar}}}
\def\LV{{\sf{LV}}}
\def\LVF{{\sf{LVF}}}
\def\BP{{\sf{BP}}}
\def\BPFP{{\sf{BPFP}}}
%\def\appleq{{<\atop{\sim}}}

% Para usar en dibujos con Pictex %%
\def\bolita{\scriptscriptstyle\bullet}

%%%%%%%
% Some abreviations for the bib database
%%%%%%%

\def\JACM{Journal of the ACM}
\def\CACM{Communications of the ACM}
\def\ICALP{International Colloquium on Automata, Languages
            and Programming}
\def\STOC{annual ACM Symp. on the Theory
          of Computing}
\def\FOCS{annual IEEE Symp. on Foundations of Computer Science}
\def\SIAM{SIAM J. Comp.}
\def\SIOPT{SIAM J. Optim.}
\def\BSMF{Bulletin de la Soci\'et\'e Ma\-th\'e\-ma\-tique de France}
\def\CRAS{C. R. Acad. Sci. Paris}
\def\IPL{Information Processing Letters}
\def\TCS{Theoret. Comp. Sci.}
\def\BAMS{Bulletin of the Amer. Math. Soc.}
\def\TAMS{Transactions of the Amer. Math. Soc.}
\def\PAMS{Proceedings of the Amer. Math. Soc.}
\def\JAMS{Journal of the Amer. Math. Soc.}
\def\LNM{Lect. Notes in Math.}
\def\LNCS{Lect. Notes in Comp. Sci.}
\def\JSL{Journal for Symbolic Logic}
\def\JSC{Journal of Symbolic Computation}
\def\JCSS{J. Comput. System Sci.}
\def\JoC{J. Compl.}
\def\MP{Math. Program.}
\sloppy

\bibliographystyle{plain}

%new macros
%%%%%%%%%%%%%%%%%%%

\def\CPRi{{\rm \#P}_{\kern-2pt\R}}

\def\CC{{\mathcal C}}
\def\DD{{\mathcal D}}
\def\NN{{\mathcal N}}
\def\MM{{\mathcal M}}
\def\GG{{\mathcal G}}
\def\UU{{\mathcal U}}
\def\ZZ{{\mathcal Z}}
\def\PP{{\mathscr P}}
\def\AA{{\mathscr A}}
\def\scC{{\mathscr C}}
\def\scE{{\mathscr E}}
\def\scD{{\mathscr D}}
\def\mcB{{\mathscr B}}
\def\mcT{{\mathscr T}}
\def\mcW{{\mathscr W}}
\def\msU{{\mathscr U}}
\def\sG{{\mathscr G}}
\def\mZ{{\mathcal Z}}
\def\mI{{\mathcal I}}
\def\sH{{\mathscr H}}
\def\sF{{\mathscr F}}
\def\bD{{\mathbf D}}
\def\nrr{{\#_{\R}}}
\def\Sd{{\Sigma_d}}
\def\oPp{\overline{P'_a}}
\def\oP{\overline{P_a}}
\def\oDH{\overline{DH_a^\dagger}}
\def\oH{\overline{H_a}}
\def\ua{\overline{\a}}
\def\Hd{\HH_{(\mathbf d)}}
\def\Hu{\HH_{\bu}}
\def\Hdi{{H_{d_i}}}
\def\Lg{{\rm Lg}}
\def\cU{\mathcal{U}}
\def\Vp{V_{\proj}}
\def\rk{\mathsf{rank}}
\def\cM{{\mathscr M}}

\def\mun{\mu_{\rm norm}}
\def\mum{\mu_{\max}}

%%%%%

\def\bx{{\bf x}}
\def\ii{{\'{\i}}}

\def\P{\mathbb P}

\newcommand{\binomial}[2]{\ensuremath{{\left(
\begin{array}{c} #1 \\ #2 \end{array} \right)}}}

\newcommand{\HH}{\ensuremath{\mathcal H}}
\newcommand{\diag}{\mathbf{diag}}
\newcommand{\CH}{\mathsf{CH}}
\newcommand{\Cone}{\mathsf{Cone}}
\newcommand{\SCH}{\mathsf{SCH}}

\newcounter{line}
\newcounter{algorithm}

\newenvironment{algorithm}[3]
{
\addtocounter{algorithm}{1}
{\bf Algorithm \thealgorithm : \sf #1} \\
{\bf Input: }#2 \\
{\bf Output: }#3 \\
\begin{list}{\arabic{line}:}{\usecounter{line}}\setlength{\leftmargin}{3em}
}
{
\end{list}
}

\newcommand{\macheps}{\varepsilon_{\mathrm{m}}}
\newcommand{\sgn}{\mathrm{sgn}}

\definecolor{orange}{rgb}{1,0.5,0}
\newcommand{\Ve}[1]{\textcolor{orange}{#1}}

\date{}

\begin{title}
{\LARGE {\bf Fast Computation of Zeros of Polynomial
Systems with Bounded Degree under Finite-precision}}
\end{title}

\author{Ir\'{e}n\'{e}e Briquel\\
Dept. of Mathematics\\
City University of Hong Kong\\
HONG KONG\\
e-mail: {\tt  irenee.briquel@gmail.com}
\and
Felipe Cucker\\
Dept. of Mathematics\\
City University of Hong Kong\\
HONG KONG\\
e-mail: {\tt macucker@cityu.edu.hk}
\and
Javier Pe\~na\\
Carnegie Mellon University\\
Tepper School of Business\\
PA, USA\\
e-mail: {\tt  jfp@andrew.cmu.edu}
\and
Vera Roshchina\\
Centro de Matem\'{a}tica e Aplica\c{c}\~{o}es\\
Universidade de \'{E}vora\\
PORTUGAL\\
e-mail: {\tt vera.roshchina@gmail.com}
}

%\date{}
\makeatletter
\maketitle
\makeatother

\begin{quote}
{\small
{\bf Abstract.}
A solution for Smale's 17th problem, for the case of
systems with bounded degree was recently given. This
solution, an algorithm computing approximate zeros
of complex polynomial systems in average polynomial time,
assumed infinite precision. In this paper we describe
a finite-precision version of this algorithm. Our main
result shows that this version works within the same
time bounds and requires a precision which, on the average,
amounts to a polynomial amount of bits in the mantissa
of the intervening floating-point numbers.
}
\end{quote}

\begin{quote}
\end{quote}

\begin{quote}
\end{quote}

%\tableofcontents

\section{Introduction}

The 17th of the problems for the 21st century
posed by Steve Smale~\cite{smale:00} asks for an algorithm
computing an approximate zero of a polynomial system
in average polynomial time.

The problem had occupied Smale during the 1990's and
led him, together with Mike Shub, to a series of
papers~\cite{Bez1,Bez2,Bez3,Bez4,Bez5,Bez6}
---known as {\em the B\'ezout series}---
where a number of ideas and results approaching a
solution for the 17th problem were proposed. These
ideas are at the core of all further research done on
Smale's problem.

Paramount within this research is the work of
Carlos Beltr\'an and Luis Miguel
Pardo~\cite{BeltranPardo08,BePa08a,BePa08b}
who provided
a randomized algorithm computing the desired
approximate zero in average expected polynomial time.
Here the word ``average'' refers to expectation over
the input data and the word ``expected'' to expectation
over the random choices made by the
algorithm\footnote{Although technically similar,
there is a remarkable
difference between the probability distributions considered.
The one for the input data is, explicitly or otherwise,
claiming some closeness to the distribution of data ``in practice.''
The only requirement for the distribution for
the random choices of the algorithm is, in
contrast, that it will be efficiently computable.
An undisputed merit of the work of Beltr\'an and Pardo
is to come up with one distribution which is so and, at the
same time, allows one to derive complexity bounds.}.
One can say that they gave a probabilistic solution to
Smale's 17th problem. Further results, including a deterministic
algorithm working in average time $N^{\Oh(\log\log N)}$
---referred to as ``nearly polynomial''--- are given
in~\cite{BC09}. This deterministic algorithm, when restricted
to systems with bounded (or even moderately growing)
degree, becomes an average polynomial-time algorithm,
refered to in~\cite{BC09} as \MD.

All the mentioned work (as well as all the work on
Smale's 17th problem not mentioned above) assume
infinite precision. As Beltr\'an and Pardo put it
in~\cite[p.~6]{BePa08b}
\begin{quote}
With the assumption of exact arithmetic [\dots] the homotopy
method [\dots] is guaranteed to produce an approximate zero
of $f$.
\end{quote}
This statement begs the question, what can one do if (as
it happens with digital computers) only finite
precision is available\footnote{Incidentally, finite precision analysis for
algorithms dealing with multivariate polynomial systems was
pioneered by Steve Smale, dragging on the way one of the
authors of the present paper, in~\cite{CS98}.}?
The goal of the present paper is to give an answer
for systems of moderate degree.

The distinctive feature of a finite precision algorithm is the
presence of a real number $u\in(0,1)$, called {\em round-off unit},
with the property that all
occurring numbers $x$ in the computation are replaced
by a number $r_u(x)$ (the rounding of $x$) such that
$|x-r_u(x)|\leq u|x|$ (a more detailed account on finite-precision
computations is in~\S\ref{subsec:bf}).
Algorithms where $u$ remains fixed
through the computation are said to have {\em fixed precision}.
Otherwise, they have {\em variable precision}. In this paper, we describe
and analyze finite-precision versions for both settings
---we denote them by \MDFix\ and \MDF, respectively---
of algorithm \MD. The rationale for the
consideration of both settings will be made clear soon enough.

Our results are of a probabilistic nature (often refered to as an
``average-case analysis''). They therefore require a probability
measure on the space of data. We next describe the measures
we use.

Let $\Hd$ denote the linear space of complex polynomial systems
$f=(f_1,\ldots,f_n)$
with $f_i$ homogeneous of degree $d_i$ in $n+1$ variables.
Given a round-off unit $\bu$, each system $f\in\Hd$ is rounded to
a system $r_{\bu}(f)$. The data space corresponding to this precision
is therefore
$$
   \Hu:=\{r_{\bu}(f)\mid f\in\Hd\}.
$$
Clearly we have a surjective function $r_{\bu}:\Hd\to\Hu$.
The measure we endow $\Hu$ with is the push-out $\nu_{\bu}$ of the standard
Gaussian $\mu$ in $\Hd$. That is, we endow $\Hd$ with a standard Gaussian
distribution (with respect to the Bombieri-Weyl basis in $\Hd$,
see \S\ref{subsec:notations} for details) and we define, for all Borelian
subset $A\subseteq \Hu$, $\nu_{\bu}(A):=\mu(r^{-1}(A))$.
With the measure $\nu_{\bu}$ at hand we can state our first main result.
\medskip

\noindent
{\bf Theorem~A.}\quad{\em
There exists a fixed precision algorithm \MDFix\ satisfying the
following.
For a random input $f\in\Hu$, algorithm \MDFix\ returns an approximate
zero of $f$ with probability at least
$$
    1-\Oh\left(\frac{D^3N(n+1)^{D+1}}{\log(1/\bu)}\right).
$$
Otherwise, \MDFix\ returns a failure message. The %expectation of the
number of arithmetic operations performed is bounded as
$$
\Oh\left(\frac{1}{n \sqrt{D}(\log N + D + n^2)\bu} \right).
$$
%$\Oh\big(D^3N^2(n+1)^D\big)$. The worst-case is bounded by an
%exponential function of $N$.

Here $N=\dim_{\C}\Hd$ denotes the size of the input systems
and $D:=\max\{d_1,\ldots,d_n\}$.
}
\bigskip

While the consideration of fixed precision is a realistic approach, a
result such as Theorem~A is not without shortcomings. Note
that \MDFix\ does not always return an approximate
zero of its input $f$, and it fails to do so with positive
probability. And, unfortunately, the lower bound for the
probability of success shown in Theorem~A becomes meaningless
when $N$ grows. In other words, a fixed precision $\bu$ puts limits
on the size of the systems for which one can expect \MDFix\ to
succeed.

The relationship between input size, round-off unit, and probability
of success implicit in Theorem~A can be better expressed with the use
of variable precision and the possibility of continuously adjusting
the input reading to the current round-off unit. This is a less
realistic setting but it pays off in terms of understanding. We
next see how.

We assume variable precision. That is, algorithms now start with
an initial round-off unit and they have the capability to refine
this parameter as the execution proceeds. In this context we will be
interested in the smallest value $u_*$ attained by the round-off
unit during the execution. A bound on $u_*$ amounts to a bound
on the number of bits (or digits)
required to store floating point approximations of the
complex numbers occurring during the computation and,
in this sense, is related to the bit-cost of performing
the computation. In fact, the maximum number of bits we will need
for each such approximation is essentially
$|\log u_*|$.

An issue naturally raised by the assumption of variable precision
is the space of data from which algorithms will take their inputs.
The spaces $\Hu$ are appropriate for the fixed precision setting:
systems $f$ are read with the precision used throughout the
computation. They do not appear to be so for the variable
precision context. But finite-precision algorithms cannot take
inputs with infinite-precision entries. An elegant solution for this
situation is the consideration of black-boxes. These are theoretical
devices which, as we said, are not realistic. But they do allow for
the statement of results highlighting in a clear manner the
relationship between precision needed and input size.

In what follows, to every $f\in\Hd$ we associate a routine
{\tt read\_input$_f$} such that
{\tt read\_input$_f$(\ )} returns an approximation of $f$
with the current round-off unit $u$. It is this routine what is given
as input to our variable precision algorithm \MDF. Because of the bijection
between this set of routines and the space $\Hd$ we may (and will)
abuse language and take $\Hd$ as the space of data for \MDF.
In particular, we will endow this space of data with the standard
Gaussian mentioned above, which we will denote by $N(0,\Id)$.

Our second main result is the following.
\medskip

\noindent
{\bf Theorem~B.}\quad{\em
%\begin{theorem}\label{thm:main2}
%Let $N(0,\Id)$ denote the standard Gaussian in the
%space $\Hd$ of polynomial systems $f=(f_1,\ldots,f_n)$
%with $f_i$ homogeneous of degree $d_i$ in $n+1$ variables.
There exists a variable precision algorithm \MDF\ satisfying the
following.
When $f$ is randomly chosen from $N(0,\Id)$,
algorithm \MDF\ on input $f$ stops almost surely,
and when it does so, returns an approximate zero of $f$.
The number of arithmetic operations
$\cost_{\MDF}(f)$ of \MDF\ on input $f$
is bounded on the average as
$$
  \E_{f\sim N(0,\Id)} \cost_{\MDF} (f) = \Oh (D^3 N^2 (n+1)^{D+1}).
$$
Furthermore, the finest precision $u_*(f)$ used by
\MDF\ on input $f$ is bounded on the average as
$$
\E_{f\sim N(0,\Id)} \log |u_*(f)| = \Oh (D^3 N (n+1)^{D+1})
$$
and as a consequence, when $D$ is bounded, the
bit-cost of \MDF\ is, on the average,
polynomial in the size $N$ of the input.
}
%\end{theorem}
\bigskip

Before proceeding with the technicalities, a couple of remarks
are in place.

{\bf (1)} As mentioned above, we confirm that Theorem~B adds
understanding to the situation depicted in Theorem~A. Indeed,
in the variable precision context the algorithm returns an approximate
zero almost surely (that is, the probability of failure is now zero).
Furthermore, the complexity
(understood as number of arithmetic operations performed)
of \MDF\ remains essentially the same as that of \MD. Finally,
the relationship between precision and input size is made clear:
we exhibit polynomial bounds (in the input size $N$) for
the expected (over random input systems $f$) number of bits
necessary to carry out the computation.

{\bf (2)} Both Theorems A and B show only the existence of
algorithms \MDFix\ and \MDF, respectively. We do not fully
exhibit these algorithms in this paper. For such a complete
description we should provide the exact values of some constants
occurring within the `big Oh' notation in a few intermediate results
in our development. This is certainly possible but we believe that
doing so would only degrade our exposition. Because, on the one hand,
it would bring undue focus on a marginal issue and, on the other
hand, the length of the exposition would unavoidably increase in a
non-trivial way. We follow in this sense a well established tradition
in finite-precision analyses, by which the goal is the understanding
of the magnitude of the precision needed by an algorithm. In
particular, we do not intend algorithms \MDFix\ and \MDF\ to
be actually implemented. They are simply vehicles to understand the
behavior of the (already implemented) algorithm \MD\ and in this
sense we make ours the words of Wilkinson quoted by Higham to
open the tenth chapter of~\cite{Higham96}:
\begin{quote}
{\small
All too often, too much attention is paid to the precise error bound
that has been established. The main purpose of such an analysis
is either to establish the essential numerical stability of an algorithm
or to show why it is unstable and in doing so to expose what sort of
change is necessary to make it stable. The precise error bound is not
of great importance.
}
\end{quote}

\section{Preliminaries}

\subsection{Setting and Notation}\label{subsec:notations}

For $d\in\N$ we denote by $H_d$ the subspace of
$\C[X_0,\ldots,X_n]$
of homogeneous polynomials of degree $d$. For $f\in H_d$
we write
$$
   f(X) = \sum_{\alpha} {d \choose \alpha}^{1/2}\,
   a_\alpha X^\alpha
$$
where $\alpha=(\alpha_0, \dots, \alpha_n)$ is assumed to range over
all multi-indices such that $|\alpha| = \sum_{k=0}^n \alpha_k = d$,
${d\choose \alpha}$ denotes the multinomial coefficient, and
$X^\alpha:= X_0^{\alpha_0}X_1^{\alpha_1}\cdots
X_n^{\alpha_n}$.
That is, we take for basis of the linear space $H_d$ the
{\em Bombieri-Weyl} basis consisting of the monomials
${d \choose \alpha}^{1/2}X^\alpha$. A reason to do so is
that the Hermitian inner product associated to this basis is unitarily
invariant. That is, if $g\in H_d$ is given by
$g(x) = \sum_{\alpha} {d\choose \alpha}^{1/2}b_\alpha X^\alpha$,
then the canonical Hermitian inner product
$$
    \langle f,g\rangle =
   \sum_{|\alpha|=d} a_\alpha\, \overline{b_\alpha}
$$
satisfies, for all elements $\nu$ in the unitary group $\cU(n+1)$, that
$$
    \langle f,g\rangle =\langle f\circ \nu,g\circ\nu\rangle.
$$
Fix $d_1,\ldots,d_n\in\N\setminus\{0\}$ and
let $\Hd=H_{d_1}\times \ldots\times H_{d_n}$ be the
vector space of polynomial systems $f=(f_1,\ldots,f_n)$ with
$f_i\in\C[X_0,\ldots,X_n]$ homogeneous of degree $d_i$.
The space $\Hd$ is naturally endowed with a Hermitian inner product
$\langle f,g \rangle = \sum_{i=1}^n \langle f_i, g_i \rangle$.
We denote by $\|f\|$ the corresponding norm of $f\in\Hd$.

We let $N:=\dim_\C\Hd$, $D:=\max_i d_i$, and
$\Dn:=\prod_i d_i$. Also, in the rest of this paper, we
assume $d_i\geq 2$ for all $i\leq n$ (linear equations
can be easily eliminated). In particular, $D\geq 2$.

Let $\proj^n:=\proj(\C^{n+1})$ denote the complex projective space
associated to $\C^{n+1}$ and $S(\Hd)$ the unit sphere of $\Hd$.
These are smooth manifolds that naturally carry the structure of a
Riemannian manifold
(for $\proj^n$ the metric is called Fubini-Study metric).
We will denote by $\dpr$ and $\dsp$ their Riemannian distances which,
in both cases, amount to the angle between the arguments.
Specifically, for $x,y\in\proj^n$ one has
\begin{equation}\label{eq:projdist}
 \cos\dpr(x,y) = \frac{|\langle x,y\rangle|}{\|x\|\, \|y\|}.
\end{equation}
Occasionally, for
$f,g\in\Hd\setminus\{0\}$, we will abuse language and write
$\dsp(f,g)$ to denote this angle,
that is, the distance $\dsp\big(\frac{f}{\|f\|},\frac{g}{\|g\|}\big) = \dsp\big(f,g\big)$.
We define the {\em solution variety} to be
$$
  \Vp:= \{(f,\z)\in \Hd\times \proj^n \mid f\neq0 \mbox{ and } f(\z)= 0\}.
$$
This is a smooth submanifold of $\Hd\times \proj^n$
and hence also carries a Riemannian structure.
We denote by $\Vp(f)$ the zero set of $f\in\Hd$ in $\proj^n$.

By B\'ezout's Theorem, $\Vp(f)$ contains $\Dn$ points for almost all $f$.
Let $Df(\zeta)_{|T_\zeta}$ denote the restriction of the derivative of
$f\colon\C^{n+1}\to\C^n$ at~$\zeta$ to the tangent space
$T_\zeta :=\{v\in\C^{n+1}\mid \langle v,\zeta\rangle = 0\}$
of $\proj^n$ at $\zeta$.
The {\em subvariety of ill-posed pairs} is defined as
\begin{equation*}\label{eq:Sigma'}
 \Sigma'_\proj := \{(f,\z)\in \Vp\mid \rk\,Df(\zeta)_{|T_\zeta} < n\} .
\end{equation*}
Note that $(f,\z)\not\in\Sigma'_\proj$ means that $\z$ is a
simple zero of~$f$.
In this case, by the implicit function theorem, the projection
$\Vp\to\Hd,(g,x)\mapsto g$
can be locally inverted around $(f,\z)$.
The image $\Sigma$ of $\Sigma'_\proj$ under the
projection $\Vp\to\Hd$ is called the
{\em discriminant variety}.

\subsection{Approximate Zeros,
Complexity and Data Distribution}\label{se:Newton}

In~\cite{Shub93b}, Mike Shub introduced the following projective
version of Newton's method.
We associate to $f\in\Hd$
(with $Df(x)$ of rank~$n$ for some~$x$) a map
$N_f:\C^{n+1}\setminus\{0\}\to\C^{n+1}\setminus\{0\}$ defined
(almost everywhere) by
$$
       N_f(x)=x-Df(x)_{|T_x}^{-1}f(x).
$$
Note that $N_f(x)$ is homogeneous of degree~0 in~$f$ and
of degree~1 in~$x$ so that
$N_f$ induces a rational map from $\proj^n$ to $\proj^n$
(which we will still denote by $N_f$) and this map
is invariant under multiplication of $f$ by constants.

We note that $N_f(x)$ can be computed from $f$ and $x$ very
efficiently: since the Jacobian $Df(x)$ can be evaluated with $\Oh(N)$
arithmetic operations~\cite{bast:83}, one can do with a total of
$\Oh(N+n^3)=\Oh(N)$ arithmetic operations, the equality
since $d_i\geq 2$ implies $N=\Omega(n^3)$.

It is well-known that when $x$ is sufficiently close to a simple zero
$\z$ of $f$, the sequence of Newton iterates beginning at $x$ will
converge quadratically fast to $\z$. This property led Steve Smale to
define the following intrinsic notion of approximate zero.

\begin{definition}\label{def:app-zero}
By an {\em approximate zero} of $f\in\Hd$ associated with a
zero $\z\in\proj^n$ of $f$
we understand a point $x\in\proj^n$ such that the
sequence of Newton iterates
(adapted to projective space)
$$
   x_{i+1}:= N_{f}(x_{i})
$$
with initial point $x_0:=x$ converges immediately
quadratically to $\z$, i.e.,
$$
 \dpr(x_i,\z)\le \Big(\frac12\Big)^{2^i -1}\ \dpr(x_0,\z)
$$
for all $i\in\N$.
\end{definition}

It is this notion of approximation that is referred to in the
statement of Smale's 17th problem.

The last notion necessary to formally state Smale's problem
is that of `average cost'. For the cost of a computation Smale
proposes the number of arithmetic operations (this includes
comparisons and possibly square roots) performed during the
computation. In the case of a finite-precision algorithm
one needs to multiply this number by the largest number
of bits (or digits) necessary to approximate the complex numbers
occurring during the computation.%% IB dire que l'on fait aussi le boulot de simuler les racines carrées

The word `average' refers to the standard normal distribution
for the data (input) system $f\in\Hd$.  Recall, we
express an element $f\in\Hd$ as a linear combination
of the monomials in the Bombieri-Weyl basis. The
standard normal distribution corresponds to choosing
the coefficients in this combination independently and
identically distributed from the centered Gaussian distribution
on $\C$ (which in turn amounts to draw real and imaginary
parts independently from the centered Gaussian distribution
on $\R$). We denote this distribution on $\Hd$ by
$N(0,\Id)$.

Hence, if $\cost(f)$ denotes the cost of computing an
approximate zero for $f$ with a given algorithm
then the average cost of this algorithm, for
inputs in $\Hd$, is given by the expected value
$$
    \E_{f\sim N(0,\Id)}\cost(f).
$$
We remark that if the cost is homogeneous of degree zero, that is,
if $\cost(f)=\cost(\lambda f)$ for all $\lambda\neq0$, then the
expectation above is the same as the expectation with $f$ drawn
from the uniform distribution on the unit sphere $S(\Hd)$.

Smale's 17th problem asks for an algorithm computing an
approximate zero (in the sense of Definition~\ref{def:app-zero})
with average cost (for the cost and data distribution
described above) bounded by $N^{\Oh(1)}$.

\subsection{Condition Numbers}\label{se:cond-num}

How close need $x$ to be from $\z$ to be an approximate zero? This
depends on how well conditioned the zero $\z$ is.

For $f\in\Hd$ and $x\in\C^{n+1}\setminus\{0\}$ we define the {\em
(normalized) condition number}~$\mun(f,x)$ by
$$
  \mun(f,x):= \|f\|\left\| \left(Df(x)_{|T_x}\right)^{-1}
  \diag\left(\sqrt{d_1}\|x\|^{d_1-1},\dots,\sqrt{d_n}
   \|x\|^{d_n-1}\right)\right\|,
$$
where the right-hand side norm denotes the spectral norm
and $\diag(a_i)$ denotes the diagonal matrix with entries~$a_i$.
Note that $\mun(f,x)$ is homogeneous of degree 0 in both arguments,
hence it is well defined for $(f,x)\in S(\Hd)\times \proj^n$. Also,
it is well known (see~\cite[Ch.~12, Corollary~3]{bcss:95})
that $\mun(f,x)\geq 1$.

The following result (essentially, a $\gamma$-Theorem in Smale's
theory of estimates for Newton's method~\cite{Smale86}) quantifies
our claim above (see \cite{BC09} for its proof).

\begin{theorem}\label{thm:gamma}
Assume $f(\zeta)=0$ and $\dpr(x,\zeta)\leq
\frac{\nu_0}{D^{3/2}\mun(f,\zeta)}$ where
$\nu_0:=3-\sqrt{7}\approx 0.3542$.
Then $x$ is an approximate zero of $f$ associated
with $\z$.\eproof
\end{theorem}

The next result, Proposition 4.1 from~\cite{BC09},
gives bounds on the variation of
the condition number $\mun(f,x)$ when $f$ and $x$ vary.

\begin{proposition}\label{prop:apps}
Assume $D\geq 2$. Let $0<\e\leq 0.13$ be arbitrary and
$C\leq \frac{\e}{5.2}$.
For all $f,g\in S(\Hd)$
and all $x,y \in \C^{n+1}$, if
$\dsp(f,g) \leq \frac{C}{D^{1/2}\mun(f,x)}$ and $\dpr(x,y) \leq
 \frac{C}{D^{3/2}\mun(f,x)}$,
then
\begin{equation}\tag*{\qed}
   \frac{1}{1+\e}\, \mun(g,y)\leq \mun(f,x) \leq
     (1+\e)\mun(g,y).
\end{equation}
\end{proposition}

In what follows, we will fix the constants
$\e:=0.13$ and $C:=\frac{\e}{5.2} = 0.025$.

We also introduce the \emph{mean square condition number} of $q$
given by
\begin{equation}\label{eq:mumean}
\mu_2^2 (q) := \frac{1}{\Dn} \sum_{\z : q(\z)=0} \mun^2 (q,\z).
\end{equation}

\subsection{An Adaptive Homotopy Continuation}

Suppose that we are given an input system $f\in S(\Hd)$ and a
pair $(g,\z)\in\Vp$, where $g$ is also in the unit sphere and such
that
$f$ and $g$ are $\R$-linearly independent.
Let $\a=\dsp(f,g)$.
Remark that one can compute $\a$ as
\begin{equation}\label{eq:alpha}
\a= 2 \arcsin \left(\frac{\|f-g\|}{2}\right).
\end{equation}

Consider the line segment
$E_{g,f}$ in $\Hd$ with endpoints~$g$ and~$f$.
We parameterize this segment by writing
$$
  E_{g,f}=\{q_\tau\in\Hd\mid \tau\in[0,1]\}
$$
with $q_\tau$ being the only point in $E_{g,f}$ such that
$\dsp(g,q_\tau)=\tau \a$.
Explicitly, as remarked in~\cite{BC09}, we have
$q_\tau = tf+(1-t)g$, where $t=t(\tau)$ is given by
\begin{equation}\label{eq:lambda}
    t(\tau) =\frac{1}{\sin\a\cot(\tau\a)-\cos\a+1}.
\end{equation}
If $E_{g,f}\cap\Sigma=\emptyset$, and hence almost surely,
this segment can be lifted to a path
given by a continuous function $[0,1]\to\Vp$ mapping
$\tau\mapsto (q_\tau,\z_{\tau})$.

In order to find an approximation
of the zero $\z_1$ of $f=q_1$ we may start with the zero $\z=\z_0$
of $g=q_0$ and numerically follow the path $(q_\tau,\z_\tau)$ by
subdividing $[0,1]$ into points $0=\tau_0<\tau_1<\cdots<\tau_k=1$
and by successively computing approximations~$x_i$ of~$\z_{\tau_i}$
by Newton's method.

This course of action is the one proposed in the B\'ezout series
and further adopted
in~\cite{BeltranPardo08,BePa08a,BePa08b,BC09}.
The (infinite precision) continuation procedure
is the following (here $\lambda=\frac{C(1-\e)}{2(1+\e)^4}
\approx 6.67\cdot 10^{-3}$,
see~\cite{BC09}).
\begin{center}
\algo
\> Algorithm \ALH\\[2pt]
\>{\bf input} $f$, $g$, $\z$\\[2pt]
\> \#\# $(g,\z)\in V$, $f\neq g$ \#\# \\ [2pt]
\>\>$\alpha:=\dsp(f,g)$, $\tau:=0$, $q_\tau:=g$\\[2pt]
\>\> {\tt repeat}\\[2pt]
\>\>\>$\Delta\tau:= \frac{\lambda}
 {\alpha D^{3/2}\mun^2(q_\tau,x)}$\\[2pt]
\>\>\>$\tau:=\min\{1,\tau+\Delta\tau\}$\\[2pt]
%\>\>\>$t:=\frac{s}{r\sin\a\cot(\tau\alpha) - r\cos\alpha + s}$\\[2pt]
\>\>\>$q_\tau:=t(\tau) f+(1-t(\tau)) g$\\[2pt]
\>\>\>$x:=N_{\tilde q_\tau}(x)$\\[2pt]
\>\>\>$x: = x/\|x\|$\\[2pt]
\>\> {\tt until} $\tau= 1$\\[2pt]
\>\> {\tt RETURN $x$}
\falgo
\end{center}
Note that the step-length $\Delta\tau$ depends on
$\mun(q_\tau,x)$. Hence, the adaptiveness.

The algorithm \MD\ (Moderate Degree) from~\cite{BC09} is
a direct application of \ALH\ having as initial pair $(g,\z)$
the pair $(\oU,\bz_1)$, where
$\oU=(\oU_1,\ldots,\oU_n)\in S(\Hd)$ with
$\oU_i = \frac1{\sqrt{2n}}(X_0^{d_i} - X_i^{d_i})$
and $\bz_1=\frac1{\sqrt{n+1}}(1,\ldots,1)$.

\algo
\> Algorithm \MD\\[2pt]
\>{\bf input} $f\in \Hd$\\[2pt]
\>\> {\tt run \ALH\ on input $(f,\oU,\bz_1)$} \falgo

\subsection{Roadmap}

Theorems~A and~B are proved by designing finite-precision
versions of algorithm \ALH\ which take into
account the errors due to the use of finite-precision.
The variable precision version \ALHF\ is
described in detail in Section~\ref{sec:homotopy}. In particular,
its main properties are shown in Theorem~\ref{thm:main1} in this section.
Proposition~\ref{prop:Recurence}
---a finite precision version of the
inductive proof of \cite[Theorem 3.1]{BC09}---
provides the backbone for the proof
of Theorem~\ref{thm:main1}.

Once with Theorem~\ref{thm:main1} at hand, the proof of Theorem~B is
carried out in a more or less straightforward manner in
Section~\ref{sec:THB}.

The use of fixed precision poses less demands in
algorithmic design (the issue of round-off unit updating
now becoming irrelevant). Algorithm \MDFix\ is therefore
a simplification of \MDF, which we describe in Section~\ref{sec:THA}
together with the proof of Theorem~A.

As just mentioned, the backbone of all this development
is Proposition~\ref{prop:Recurence}. The proof of this result
relies on finite-precision estimates for the errors in a number of
basic procedures. These estimates are collected in the
next section. We are aware they do not make the most exciting
part of the paper but it is a part we cannot do without.

\section{Error Bounds}\label{sec:errorbounds}

In this section we show bounds for the basic computations
occurring in \ALHF. We will
use these bounds in subsequent sections to show our main
result.

\subsection{Basic facts}\label{subsec:bf}

We recall the basics of a floating-point arithmetic which idealizes
the usual IEEE standard arithmetic. In contrast to the standard model
(as in \cite{Higham96}) we adapt our exposition to complex arithmetic.
This system is defined by a set
$\F\subset\Q[i]$ containing $0$ (the {\em floating-point complex numbers}), a
transformation $r_u:\C\to\F$ (the {\em rounding map}), and a constant
$u\in\R$ (the {\em round-off unit}) satisfying $0<u<1$. The
properties we require for such a system are the following:
\begin{enumerate}
\item[(i)]
For any $x\in\F$, $r_u(x)=x$. In particular, $r_u(0)=0$.
\item[(ii)]
For any $x\in\C$, $r_u(x)=x(1+\d)$ with $|\d|\leq u$.
\item[(iii)]
For any $y\in\F$, the set $r_u^{-1}(y)$ is measurable in $\C$.
\end{enumerate}
Property~(iii) ensures that the measure $\nu_{\bu}$ described in the 
Introduction is well defined. Because of the enumerability of $\Q[i]$, 
this measure can be seen as a discretization of the Gaussian in $\Hd$. 
\smallskip

We also define on $\F$ arithmetic operations following the classical
scheme
$$
  x\tilde\circ y=r_u(x\circ y)
$$
for any $x,y\in\F$ and $\circ\in\{+,-,\times,/\}$, so that
$$
  \tilde\circ:\F\x\F\to\F.
$$

The following is an immediate consequence of property (ii) above.

\begin{proposition}\label{prop:f1}
For any $x,y\in\F$ we have
\begin{equation}\tag*{\qed}
        x\tilde\circ y=(x\circ y)(1+\d),\qquad |\d|\leq u.
\end{equation}
\end{proposition}

When combining many operations in floating-point arithmetic,
quantities such as $\prod_{i=1}^n(1+\d_i)^{\rho_i}$ naturally
appear.
Our round-off analysis uses the notations and ideas
in Chapter~3 of~\cite{Higham96}, from where we quote
the following results:

%The proof of the following propositions can be found in
%Chapter~3 of~\cite{Higham96}.
%The notation they introduce, the
%quantities $\g_n$ and $\t_n$, and the relations showed therein, will
%be used in our round-off analysis.

\begin{proposition}\label{propHigham}
If $|\d_i|\leq u$, $\rho_i\in\{-1,1\}$, and $nu<1$, then
$$
  \prod_{i=1}^n(1+\d_i)^{\rho_i}=1+\theta_n,
$$
where
\begin{equation}\tag*{\qed}
        |\theta_n|\leq \g_n=\frac{nu}{1-nu}.
\end{equation}
\end{proposition}

\begin{proposition}\label{propHigham2}
For any positive integer $k$ such that $ku<1$, let $\theta_k, \theta_j$ be any
quantities satisfying
$$
        |\theta_k|\leq \g_k=\frac{ku}{1-ku}
\hspace{3em}
        |\theta_j|\leq \g_j=\frac{ju}{1-ju}
    .
$$
The following relations hold.
\begin{enumerate}
\item[{\rm 1.}]
$(1+\theta_k)(1+\theta_j)=1+\theta_{k+j}$ for some $|\theta_{k+j}|
\le \g_{k+j}$.

\item[{\rm 2.}]
$$
\frac{1+\theta_k}{1+\theta_j}=\left\{
\begin{array}{ll}
        1+\theta_{k+j}& \mbox{if $j\leq k$,}\\
        1+\theta_{k+2j}& \mbox{if $j> k$.}
\end{array}\right .
$$
for some $|\theta_{k+j}|
\le \g_{k+j}$ or some $|\theta_{k+2j}|
\le \g_{k+2j}$.
\item[{\rm 3.}]
If $ku,ju\leq1/2$, then $\g_k\g_j\leq \g_{\min\{k,j\}}$.

\item[{\rm 4.}]
$i\g_k\leq \g_{ik}$.

\item[{\rm 5.}]
$\g_k+u\leq \g_{k+1}$.

\item[{\rm 6.}]
$\g_k+\g_j+\g_k\g_j\leq \g_{k+j}$. \eproof
\end{enumerate}
\end{proposition}

From now on, whenever we write an expression containing $\theta_k$ we mean
that the same expression is true for some $\theta_k$, with
$|\theta_k| \le \g_k$.

When computing an arithmetic expression $q$ with a round-off
algorithm, errors will accumulate and we will obtain another
quantity which we will denote by $\fl(q)$.
For a complex number, we write $\Error(q)=|q-\fl(q)|$; for vectors or matrices, $\Error(q)$ will denote the vector or matrix of coordinates $|q_\alpha - \fl(q_\alpha)|$, allowing us to choose various norms to estimate this error.

An example of round-off analysis which will be useful in what
follows is given in the next proposition, the proof of which follows the lines of
the proof of the real version of this result that can be
found in Section~3.1 of~\cite{Higham96}.

\begin{proposition}\label{Cpxscalar}
There is a finite-precision
algorithm which, with input $x,y\in\C^n$,
computes the inner product of $x$ and $y$. The computed value
$\fl(\langle x,y\rangle)$ satisfies
$$
   \fl(\langle x,y\rangle)=\langle x,y\rangle+
   \t_{\lceil\log_2 n\rceil+1}\sum_{i=1}^n |x_i\bar y_i|.
$$
In particular, if $x=y$, the
algorithm computes $\fl(\|x\|^2)$ satisfying
\begin{equation}\tag*{\qed}
   \fl(\|x\|^2)=\|x\|^2(1+\t_{\lceil\log_2 n\rceil+1}).
\end{equation}
\end{proposition}

We assume that, besides the four basic operations, we are allowed to
compute basic trigonometric functions (such as $\sin$ and $\cos$) and the
square root with finite precision. That is, if $\op$ denotes any of these two operators, we compute $\widetilde{\op}$ such that
$$
\widetilde{\op}(x) = \op(x)(1+\delta), |\delta|<u.
$$

The following sensitivity results will help us to deal
with errors in computing trigonometric functions.

\begin{lemma}\label{lem:trigsensitivity}

\begin{enumerate}
\item[(i)] Let $t,\theta\in \R$. Then
$$
|\cos(t+\t)-\cos t| \leq |\t|;
$$
$$
|\sin(t+\t)-\sin t| \leq |\t|;
$$
\item[(ii)]
Given two reals $a$ and $e$ such that
both $a$ and $a+e$ are in the interval
$[0,0.8],$ one has
$$
|\arcsin (a+e) - \arcsin(a)| \leq 2 |e| \mbox{, with }|v| \leq |e|.
$$
\end{enumerate}
\end{lemma}
\proof

{\bf (i)}
Observe that
$$
\left|\cos (t+\t) - \cos t\right| = 2 \left|\sin\left(t+\frac{\t}{2}\right)\right|\left|\sin\frac{\t}{2}\right|\leq 2 \left|\sin\frac{\t}{2}\right|\leq |\t|,
$$
and analogously
\begin{equation}\tag*{\qed}
\left|\sin (t+\t) - \sin t\right|
= 2 \left|\cos\left(t+\frac{\t}{2}\right)\right|
\left|\sin\frac{\t}{2}\right|\leq |\t|.
\end{equation}

{\bf (ii)}
Without loss of generality, let us suppose that $e>0$.

From the intermediate value theorem, there exists a $\xi$ in $[a,a+e]$ such that
$\arcsin (a+e) = \arcsin(a) + e \arcsin'(\xi) = \arcsin(a) + e \frac{1}{\sqrt{1- \xi^2}}$.

Since $\xi \in [a,a+e]$, $ |\xi | \leq 0.8$ and thus $|\arcsin'(\xi)| \leq \frac{1}{\sqrt{1-0.8^2}} < 2$.
\eproof

To avoid burdening ourselves with the consideration
of multiplicative constants, we introduce a further
notation. Computational errors in our context are functions
of the integer parameters $n,N$ and $D$ as well as on the
condition $\mun(g,z)$ of the initial pair $(g,z)$. For
any such function $\Phi$, we will write
$$
\ll \Phi \rr := \theta_{\Oh( \Phi)}.
$$

The next properties follow directly from the properties of the $\theta$ notation.
\begin{proposition}\label{prop:calcbigOh}
Let $\Phi$ and $\Psi$ be two real functions. The following relations hold:
\begin{enumerate}
\item[{\rm 1.}]
$\ll \Phi \rr + \ll \Psi \rr = \ll \max (\Phi,\Psi) \rr.$
\item[{\rm 2.}]
$\ll \Phi \rr \ll \Psi \rr = \ll \max (\Phi,\Psi) \rr.$
\item[{\rm 3.}]
If $\Phi \geq 1$, $ \Phi \ll \Psi \rr = \ll \Phi \Psi \rr.$
\end{enumerate}
\end{proposition}

\subsection{Bounding errors for elementary computations}

We now begin showing bounds for the errors in
the crucial steps of our algorithm. To avoid burdening the
exposition we will do so only for
the steps dominating the accumulation of errors and
simply warn the reader of the minor steps we consider
as exact. %For instance, we suppose, that we can compute
%square roots and some trigonometric functions
%with finite precision, similarly as for the four basic
%arithmetic operations.

%For instance, in all what follows we
%assume that $\alpha=\dsp(f,g)$ is computed exactly.
%Similarly, in the next result we assume that
%$1/(2 \sin\frac{\a}{2})$ is computed exactly.

We begin with the evaluation of the errors in computing $\a$.
Remark that we suppose $\alpha \leq \pi/2$ in the following lemma.
This will be ensured by the computation of $\alpha$ at the beginning of
\ALHF. If this quantity is more than $\pi /2 $, we set $f=-f$, ensuring that
$\alpha \leq \pi /2$. We neglect the errors in this operation, and thus suppose in the remainder that
$\alpha \leq \pi /2$.

\begin{lemma}~\label{lemma:alpha}
Given $f$ and $g$ in $S(\Hd)$ such that $\dsp (f,g) \leq \pi/2$, one can compute $\alpha = \dsp (f,g)$ with finite
precision such that
$$
\fl (\alpha)= \alpha (1+\ll \log N\rr).
$$
\end{lemma}

\proof
As remarked in~\eqref{eq:alpha}, one can compute $\a = \dsp (f,g)$ as $\a= 2
\arcsin\left( \frac{\|f-g\|}{2}\right)$.

We can compute the norm $\|f-g\|$ similarly as the vector norm in Proposition~\ref{Cpxscalar}. In the case of
polynomials in $\Hd$, the sum is over $N$ coefficients, and thus we prove similarly that
$\fl (\|f-g\|^2) = \| f - g\|^2 (1 +  \t_{\lceil \log N \rceil +1} )$.
Since we supposed that we can compute square root with finite precision,
we get
$$
\fl(\| f-g \|) = \| f - g\| (1 +  \t_{\lceil \log N \rceil +2} ).
$$

Remark that, since we supposed $\dsp(f,g) \leq \pi/2$ and $\|f\|=\|g\|=1$,
we have $\| f-g \|/2 \leq \sin (\pi/4) = 1 / \sqrt{2}<0.71.$
We can suppose that $u$ is small enough such that the term $\t_{\lceil \log N \rceil + 2}$
is smaller than $0.8-0.71$, and thus such that $\fl (\| f-g \|/2)$ is also in $[0, 0.8]$.
We can thus apply Lemma~\ref{lem:trigsensitivity}, and by
supposing that we are able to compute the function $\arcsin$ with finite precision,
we conclude that we can compute
$\alpha=2 \arcsin\left( \frac{\|f-g\|}{2}\right)$ such that
\begin{eqnarray*}
\fl (\alpha) & = & \left(2 \arcsin \left(\frac{\|f-g\|}{2}\right)
        + 2\frac{\|f-g\|}{2} \t_{\Oh(\log N)}\right) (1+ \t_{\Oh(1)})\\
      & = & 2 \arcsin \left(\frac{\|f-g\|}{2}\right) (1+ \ll\log N\rr),
\end{eqnarray*}
the last line since $\left|\frac{\|f-g\|}{2}\right| \leq \left|\arcsin \left(\frac{\|f-g\|}{2}\right)\right|.$
\eproof

\begin{proposition}\label{prop:Errt}
Given $\tau \in \R_+$, $f$ and $g$ in $S(\Hd)$ such that $\dsp(f,g) \leq \pi/2$,
we can calculate $t(\tau)$
with finite precision such that
$$
     \fl(t) = t (1+\ll\log N\rr).
$$
\end{proposition}
\proof
First of all, observe that
\begin{eqnarray*}
t(\tau) = \frac{1}{\sin\a \cot(\a \tau)-\cos(\a)+1}
   & = & \frac{\sin(\tau \a)}{\sin\a \cos(\tau \a)-\cos\a \sin(\tau \a)+\sin(\tau \a)}\\
   & = & \frac{\sin(\tau \a)}{\sin(\a-\tau \a)+\sin(\tau \a)}\\
   & = & \frac{\sin(\tau \a)}{2\sin\left(\frac{(1-\tau)\a+\tau \a}{2}\right)\cos\left(\frac{(1-\tau)\a-\tau \a}{2}\right)}\\
   & = & \frac{\sin(\tau \a)}{2\sin\frac{\a}{2}
   \cos\left(\left(\frac{1}{2}-\tau\right)\a\right)}.
\end{eqnarray*}

We compute $t(\tau)$ via the last equality.
First, we compute $\alpha$ following Lemma~\ref{lemma:alpha}.
Then, we show easily using Lemma~\ref{lem:trigsensitivity}
that each term in the fraction can be computed with finite precision
up to a multiplicative factor $(1+\t_{\Oh(\log N)}).$
We conclude using Proposition~\ref{propHigham2}.
\eproof

% we have
%$$
%\fl (\sin(\tau \a)) = \left[\sin (\tau\a)+ \t_1 \tau \a\right](1+\t_1) = \sin (\tau\a)\left(1+ \frac{\t_1 \tau \a}{\sin (\tau\a)}\right)(1+\t_1) .
%$$
%Notice that since $\sin (\tau \a)\geq \tau \a -\frac{\tau^3 \a^3}{3!}$, and $\a\leq 2$, we have
%$$
%\frac{\tau \a}{\sin (\tau\a)}\leq \frac{1}{1-\frac{4\tau^2}{6}} \leq 3,
%$$
%and hence
%$$
%\fl (\sin(\tau \a)) =\sin(\tau \a)(1+\t_{4}).
%$$
%Similarly, using Lemma~\ref{lem:trigsensitivity} for the denominator we have
%\begin{eqnarray*}
%\fl\left(\cos\left(\left(\frac{1}{2}-\tau\right)\a\right)\right)
%    & = & \cos\left(\left(\frac{1}{2}-\tau\right)\a\right)\left(1+\frac{\t_2\left(\frac{1}{2}-\tau\right)\a}{\cos\left(\left(\frac{1}{2}-\tau\right)\a\right)}\right)(1+\t_1).
%\end{eqnarray*}
%Since $\left|\left(\frac{1}{2}-\tau\right)\alpha\right|\leq 1$, we have
%$$
%\frac{\left(\frac{1}{2}-\tau\right)\a}{\cos\left(\left(\frac{1}{2}-\tau\right)\a\right)}\leq \frac{1}{\cos 1}\leq 2.
%$$
%Therefore,
%$$
%\fl\left(\cos\left(\left(\frac{1}{2}-\tau\right)\a\right)\right) = \cos\left(\left(\frac{1}{2}-\tau\right)\a\right)(1+\t_{5}).
%$$
%To finish the proof, it remains to evaluate the impact of the remaining division operation, and hence
%\begin{equation}\tag*{\qed}
%\fl(t)= t(1+\t_9).
%\end{equation}

% Bound on the quantity f(x)=\sum c_J X^J when ||f||=1
The following lemma bounds by $\|q \|$ the value of a polynomial
$q$ at any point on the unit sphere.

\begin{lemma}~\label{lemma:BoundBWNorm}
Given $d \in \N$, $q \in \C[X_0,\ldots,X_n]$ homogeneous
of degree $d$ and $x\in S(\C^{n+1})$, we have $|q(x)| \leq \|q\|$.
\end{lemma}

\proof
Since our norm $\|\ \|$ on $\C[X_0,\ldots,X_n]$
is unitarily invariant, for each
element $\phi \in \cU(n+1)$, one has $\|f \circ \phi \| = \|f \|$.

Let $e_0:=(1,0,\ldots,0)\in\C^{n+1}$.
Taking $\phi$ such that $\phi(e_0) = x$, one has
$$
|q(x)| = | (q \circ \phi \circ \phi^{-1}) (x) | = | (q \circ \phi) (e_0)|.
$$

But $| q \circ \phi (e_0)|$ is exactly the coefficient of $X_0^d$
in $q\circ\phi$ with respect to the Bombieri-Weyl
basis of $\C[X_0,\ldots,X_n]$, and thus
$| q \circ \phi (e_0)| \leq \|q \circ \phi\|= \|q\|$.
\eproof

\begin{proposition}\label{prop:Errq}
Given $q\in S(\Hd)$ and $x\in S(\C^{n+1})$, we can compute
$q(x)$ with finite precision $u$ such that
$$
   \|\Error(q(x))\| = \ll \log N+D \rr.
$$
\end{proposition}

\proof
For $i\leq n$, write $q_i(x) = \sum c_J x^J$. To
compute $q_i(x)$ we compute each monomial $c_J x^J$
first, and then evaluate the sum. We have
\begin{equation*}
\fl(c_J x^J) = c_J x^J(1+\t_{d_i+1}),
\end{equation*}
and thus $\Error (c_J x^J) \leq |c_J| |x|^J \g_{d_i +1}$.

As
$$
\fl(q_i (x)) = \fl \left( \sum c_J x^J \right),
$$
using pairwise summation (see section 4.2 in \cite{Higham96}) we have
\begin{eqnarray*}
\Error(q_i(x))
    &=&  \Bigl|\sum\fl (c_J x^J)-\sum (c_J x^J) \\
    &&+\sum  \fl(c_J x^J) \t_{\lceil\log_2 N\rceil}\Bigr|\\
    &\leq & \sum\Error(c_J x^J)+\sum | c_J x^J |\g_{\lceil\log_2 N\rceil}\\
    &&+\sum \Error(c_J x^J)\g_{\lceil\log_2 N\rceil}\\
    &\leq & \sum |c_J| |x|^J (\g_{D+1} + \g_{\lceil\log_2 N\rceil} + \g_{D+1} \g_{\lceil\log_2 N\rceil})\\
    &\leq &\sum |c_J| |x|^J \g_{\lceil\log_2 N\rceil+D+1}.  \quad \mbox{ (by Proposition~\ref{propHigham2} {\rm 6.})}
\end{eqnarray*}

Note that
$\sum |c_J| |x|^J \leq \|q_i\|,$
by applying Lemma~\ref{lemma:BoundBWNorm} to the polynomial of
coefficients $|c_J|$, which has the same
norm as $q_i$, at the point $|x| \in S(\C^{n+1})$. Hence,
$$
 \Error(q_i(x)) \leq  \|q_i\| \g_{\lceil\log_2 N\rceil+D+1},
$$
and therefore,
\begin{equation*}
\|\Error(q(x))\|^2
 \leq  \g_{\lceil\log_2 N\rceil+D+1}^2 \sum_i \|q_i\|^2 = \g_{\lceil\log_2 N\rceil+D+1}^2 \|q \|^2 =\g_{\lceil\log_2 N\rceil+D+1}^2.
\end{equation*}
We finally have
\begin{equation}\tag*{\qed}
\|\Error(q(x))\|=\ll \log N+D \rr.
\end{equation}

\subsection{Bounding the error in the computation of
$\mun^{-1}(q,x)$}\label{subsection:ErrorMun}

The bounds in $\Error(\mun^{-1}(q,x))$
scale well with $q$.
Hence, to simplify
notation, in all what follows we assume $\|q\|=1$.

The main result in this subsection is the following.

\begin{proposition}\label{prop:sigmaM}
Given $q\in S(\Hd)$ and $x\in S(\C^{n+1})$
we can compute $\mun^{-1}(q,x)$ satisfying
$$
    \Error(\mun^{-1}(q,x))= \ll n (\log N + D + n)\rr.
$$
\end{proposition}

Note that under the assumption  $\|q\|=1$
our condition number becomes
$$
  \mun(q,x):= \left\| \left(Dq(x)_{|T_x}\right)^{-1}
  \diag(\sqrt{d_1},\dots,\sqrt{d_n}
   )\right\|.
$$
Given $q\in S(\Hd)$ and $x\in S(\C^{n+1})$,
let $M_q\in\C^{n\times n}$ be a matrix representing
the linear operator
\begin{equation}\label{eq:DefM}
  \left[\begin{matrix} \frac{1}{\sqrt{d_1}} \\ &
  \frac{1}{\sqrt{d_2}} \\
   & & \ddots \\ & & & \frac{1}{\sqrt{d_n}}
  \end{matrix} \right]Dq(x)_{|T_x}
\end{equation}
in some orthonormal basis of $T_x$ (note that $M_q$ depends also on $x$;
that point $x$ will always be clear from the context).
We then have
$\mun^{-1}(q,x)=\|M_q^{-1}\|^{-1}
=\sigma_{\min}(M_q)$ where $\sigma_{\min}$
denotes smallest singular value. We will compute
$\mun^{-1}(q,x)$ by computing $M_q$ and then
$\sigma_{\min}(M_q)$.

The following proposition contains several technical ideas that will
help us to deal with the matrices $D{q}(x)_{|T_x}$ and $M_q$.
We use ideas from the proof of \cite{CKMM08} modifying them to
the complex case.

\begin{proposition}\label{prop:CalcDeriv}
Let $q\in \Hd$ and $x\in S(\C^{n+1})$. Then the following
statements are true:
\begin{itemize}
  \item [(i)] The restriction of the derivative of $q$ to the tangent
   space $T_x$ can be represented by the following matrix :
  $$
     D q(x)_{|T_x} = Dq(x)H,
  $$
  where $H \in \C^{(n+1) \times n}$ is the matrix made with the
  last $n$ columns of the matrix $H_x$
  defined by
  $$
  H_x = \alpha \left(I_{n+1}- 2 yy^*\right), \qquad y =\frac{x-\alpha e_{0}}{\|x- \alpha e_{0}\|}, \qquad \alpha = \frac{x_{0}}{|x_{0}|}
  $$
  if $|x_{0}| \neq 1$, and $H_x = \alpha \,I_{n+1}$ otherwise.
  \item [(ii)]
  $$
  \left\|\diag\left(\frac{1}{\sqrt{d_1}},\dots,\frac{1}{\sqrt{d_{n}}}\right)
  Dq(x)_{|T_x}\right\| \leq \|q\|.
  $$
  \item [(iii)]
  $$
  \|Dq(x)\|_F \leq \sqrt{D} \|q\|, \qquad  \left\|Dq(x)_{|T_x}\right\|_F \leq \sqrt{D} \|q\|.
  $$
\end{itemize}
\end{proposition}

\proof
{\bf (i) }
For any unitary matrix $H_x$ such that $H_x e_{0} = x$, the $n$ last columns $H$ of $H_x$ form an
orthonormal basis of $T_x$. Thus $Dq(x)H$ is the representation of $D q(x)_{|T_x}$ in that basis.

The matrix $H_x$ computed here is constructed in~\cite{BePa08b};
one checks easily that it is
unitary and that it satisfies $H_x e_{0} = x$.

{\bf (ii) }  Let $g = q \circ H_x$. Then, differentiating the equality
$g_i(H_x^*x) = q_i(x)$ and multiplying both sides by $H$
on the right, we have
\begin{equation}\label{eq:001}
D g_i(e_0)H_x^*H = Dq_i(x) H=Dq_i(x)_{|T_x},
\end{equation}
where the last equality is by (i). Observe that $H_x^*H = [e_1,\dots,
e_n]$, hence,
\begin{equation}\label{eq:002}
  D g_i(e_{0})H_x^*H =  D g_i(e_{0})_{|T_{e_{0}}}
  = \left[\frac{\partial g_i}{\partial X_1}(e_{0}),\dots,
  \frac{\partial g_i}{\partial X_n}(e_{0})\right].
\end{equation}
If we denote
$g_i (X) = \sum_{\alpha} {d \choose \alpha}^{1/2}\, g_{i \alpha}
X^\alpha$, it is straightforward that
$$
   \frac{\partial g_i}{\partial X_j}(e_{0})
   = {d_i \choose d_i -1}^{1/2} g_{i(e_{j}+(d_i-1)e_{0})}
   = \sqrt{d_i} \cdot g_{i(e_{j}+(d_i-1)e_{0})}.
$$
Therefore, from (\ref{eq:002}),
\begin{equation}\label{eq:003}
{ \left\|\frac{1}{\sqrt{d_i}} D g_i(e_{0})_{T_{e_{0}}}\right\|}^2
= \sum_j  g^2_{i(e_{j}+(d_i-1)e_{0})} \leq  \|g_i\|^2,
\end{equation}
and hence by (\ref{eq:003}) we have
\begin{eqnarray}\label{eq:004}
{\left\| \diag\left(\frac{1}{\sqrt{d_1}},\dots,\frac{1}{\sqrt{d_{n}}}\right)
D g(e_{0})_{T_{e_{0}}}\right\|_F}^2
&\leq& \sum_{i=1}^n \left\|
   \frac{1}{\sqrt{d_i}}D g_i(e_{0})_{T_{e_{0}}}\right\|^2 \nonumber\\
&\leq& \sum \left\|g_i\right\|^2= \|g\|^2.
\end{eqnarray}

Since the Hermitian inner product associated with the Bombieri-Weyl
basis is unitarily invariant, we have
$$
      \|g\|^2=\langle g,g\rangle
      = \langle q \circ H_x,q \circ H_x\rangle = \langle q,q\rangle
      = \|q\|^2 ,
$$
which by (\ref{eq:001}),(\ref{eq:002}) and (\ref{eq:004}), and since the spectral norm of a matrix is not greater than its Frobenius norm, yields
$$
   \left\| \diag\left(\frac{1}{\sqrt{d_1}},\dots,\frac{1}{\sqrt{d_{n}}}\right)
   D q(x)_{T_{x}}\right\|= \left\| \diag\left(\frac{1}{\sqrt{d_1}},\dots,
   \frac{1}{\sqrt{d_{n}}}\right) D g(e_{0})_{T_{e_{0}}}\right\|
   \leq \|g\| = \|q\|.
$$
The relations {\bf (iii)} can be shown similarly.
\eproof

The following two statements are similar to those proved
in~\cite{CKMM08} in the real case and similar
ideas are used in the proofs.

\begin{proposition}\label{prop:ErrD}
Given $q\in S(\Hd)$ and $x\in S(\C^{n+1})$, we have
$\|Dq(x)_{|T_x}\|\leq \sqrt{D}$, and we can compute
$Dq(x)_{|T_x}$ with finite precision such that
$$
  \|\Error(Dq(x)_{|T_x})\|_F\leq
  \ll n \sqrt{D}(\log N +D +\log n)\rr.
$$
\end{proposition}

\proof
The inequality $\|Dq(x)_{|T_x}\|\leq \sqrt{D}$ follows from
$\|q\|= 1$ and Proposition~\ref{prop:CalcDeriv}(iii).

We compute $Dq(x)_{|T_x}$
as in Proposition~\ref{prop:CalcDeriv}(i).
Hence each entry $(i,j)$ of the matrix $Dq(x)_{|T_x}$
is calculated as the product of $Dq_i(x)$  and the $j$th
column $H_j = (h_{kj})_{1\leq k \leq n+1}$ of $H$.
Proceeding as in the proof of Proposition~\ref{prop:Errq}
we can compute
$\frac{\partial q_i}{\partial X_k}(x)$ with
$$
  \Error\left(\frac{\partial q_i}{\partial X_k}(x)\right)
  = \ll(\log N +d_i) \rr \| Dq_i(x)\|_F.
$$
%the factor $\sqrt{D}$ being due to the fact that
%$\| Dq(x)_{|T_x}\|\leq \sqrt{D}$.

One can compute $\alpha$  as $\frac{x_{0}}{\sqrt{ x_{0} x^*_{0}}}$ with two arithmetic operations and
one square root. Observe that to compute $x-\alpha e_{0}$, we need to perform only two more
arithmetic operations. Also,
$$
    \left(yy^*\right)_{ij} =
    \frac{1}{\|x- \alpha e_{0}\|^2}\left((x- \alpha e_{0})(x- \alpha e_{0})^*\right)_{ij}
$$
and we have
\begin{eqnarray*}
    \fl \left((x- \alpha e_{0})(x- \alpha e_{0})^*\right)_{ij}
  &=& \fl\left((x- \alpha e_{0})_i\overline{(x- \alpha e_{0})}_j\right)\\
  &= &\left((x- \alpha e_{0})(x-\alpha e_{0})^*\right)_{ij}(1+\t_{11}).
\end{eqnarray*}
Further,
\begin{eqnarray*}
   \fl \left(\|x- \alpha e_{0}\|^2\right)
    &=& \fl \left(\sum_{i=1}^{n}x_i\overline x_i
     +(x_{0}-\alpha )(\overline x_{0} -\alpha )\right)\\
    &=& \|x-e_{0}\|^2(1+\t_{\lceil\log_2{(n+1)}\rceil+ 11}).
\end{eqnarray*}
Here we used pairwise summation bounds again.

Thus, by Proposition~\ref{propHigham2}(2),
$$
    \fl \left(2 yy^*\right)_{ij} = \left(2 yy^*\right)_{ij}
    (1+\t_{2 \lceil\log_2{(n+1)}\rceil+35}).
$$

Finally, taking into account one more addition and the multiplication by $\alpha$,
we get
$$
   \Error(h_{ij}) = \theta_{2\lceil\log_2{(n+1)}\rceil+39}= \ll\log n\rr.
$$
Applying Proposition~\ref{Cpxscalar}, we conclude
\begin{eqnarray*}
\Error([Dq(x)_{|T_x}]_{ij})
    &=& |\fl(\langle D{q}_i(x), H_j \rangle) -
    \langle D{q}_i(x), H_j \rangle| \\
    &=& |\langle \fl(D{q}_i(x)), \fl(H_j) \rangle \\
    & &   + \t_{\lceil\log_2 n\rceil+1} \sum_k |
      {D{q}_i(x)}_k \overline{H_{kj}} | -
     \langle D{q}_i(x), H_j \rangle|\\
%\end{eqnarray*}
%Since the precision $u$ we use verifies (\ref{eq:BoundOnU}),
%we can neglect the terms of order greater than 1 in $u$, and write
%\begin{eqnarray*}
%\Error([Dq_\tau(x)_{|T_x}]_{ij})
    &\leq& |\langle \Error (D{q}_i(x)), H_j \rangle|
   + |\langle D{q}_i(x), \Error(H_j)\rangle| \\
    & & + |\langle \Error (D{q}_i(x)),\Error(H_j)\rangle|
    + \g_{\lceil\log_2 n\rceil+1} |D{q}_i(x)| |H_j|\\
    &=& ( \ll (\log N + D) \rr +  \ll \sqrt{n} \log n \rr \\
    & & + \ll (\log N + D) \rr \ll \sqrt{n} \log n \rr
        + \ll \log n \rr ) \|Dq_i(x)\|_F \\
    &=& \ll \sqrt{n} (\log n + \log N + D) \rr \| Dq_i(x)\|_F.
\end{eqnarray*}
%The second equality holds since $\|H_j\|\leq 1$, as $H_x$ is unitary, and $\left\|{Dq_\tau}_i(x)\right\|\leq \sqrt{D}$ by %Proposition~\ref{prop:CalcDeriv} (iii).
This implies
\begin{equation}\tag*{\qed}
   \|\Error(Dq(x)_{|T_x})\|_F =
   \ll n (\log n+D+\log N)\rr \|Dq(x)\|_F= \ll n \sqrt{D} (\log n+D+\log N)\rr
\end{equation}

\begin{proposition}\label{prop:ErrM}
Given $q\in S(\Hd)$, $x\in S(\C^{n+1})$ and
$M_q$ defined by (\ref{eq:DefM}), we have $\|M_q\|\leq 1$.
In addition, we can compute such a matrix
$M_q$ with finite precision $u$ such that
$$
    \|\Error(M_q)\|_F=\ll n (\log N + D+\log n)\rr.
$$
\end{proposition}

\proof
The inequality $\|M_q\|\leq 1$ follows directly from
Proposition~\ref{prop:CalcDeriv}(ii), as $\|q\|=1$.
Floating-point errors can be evaluated exactly as in
Proposition~\ref{prop:ErrD}; however, one gets rid of the factors
$\sqrt{D}$ since the bound on $\|M_q\|$ is better than the bound on
$\|Dq(x)_{|T_x}\|$. As a counterpart, one has to take into account
one more division by $\sqrt{d_i}$ of each entry of the matrix, which
slightly changes the constants, but leaves the order in $N,D $ and $n$
unchanged.
\eproof

\proofof{Proposition~\ref{prop:sigmaM}}
We use ideas from the proof of an analogous
proposition in~\cite{CKMM08}.
Let $x\in S(\C^{n+1})$, $q\in S(\Hd)$ and $M_q$ be as in
Proposition~\ref{prop:ErrM}. Then
$\mun^{-1}(q,x)
=\sigma_{\min}(M_q)=\|M^{-1}_q\|^{-1}$
and we can compute the first expression by
computing the last.

Let $E'=M_q-\fl(M_q)$. By Proposition~\ref{prop:ErrM},
$$
  \|E'\|\leq \|E'\|_F\leq \ll n (\log N + D +\log n)\rr.
$$
Let $\MM_q=\fl(M_q)$. We compute
$\sigma_{\min}(\MM_q)=\|M_q^{-1}\|^{-1}$ using a backward
stable algorithm (e.g., QR factorization).
Then the computed $\fl(\sigma_{\min}(\MM_q))$ is the exact
$\sigma_{\min}(\MM_q+E'')$ for a matrix $E''$ with
\[
  \|E''\| \leq cn^2 u \|\MM_q\|
\]
for some universal constant $c$ (see, e.g., \cite{GoLoan,Higham96}).
Thus,
$$
  \fl(\sigma_{\min}(M_q))=\fl(\sigma_{\min}(\MM_q))=
  \sigma_{\min}(\MM_q+E'')=\sigma_{\min}(M_q+E'+E'').
$$
Write $E=E'+E''$. Then, using $\|M_q\| \leq 1$ (by Proposition~\ref{prop:ErrM}),
\begin{eqnarray*}
  \|E\| &\leq& \|E'\|+\|E''\|
  \leq \|E'\| + cn^2 u \|\MM_q\|
  \leq \|E'\| + cn^2 u (\|M_q\|+\|E'\|)\\
  &=& \ll n (\log N + D +\log n)\rr
   + \ll n^2 \rr (1+\ll n (\log N + D +\log n)\rr)\\
   &=& \ll n (\log N + D +\log n)\rr
  + \ll n (\log N+D+ n)) \rr \\
  &=& \ll n (\log N + D + n)\rr,
\end{eqnarray*}
using Proposition~\ref{prop:calcbigOh} in the penultimate row.

Therefore, $\fl(\sigma_{\min}(M_q))
= \sigma_{\min}(M_q+E)$ which
implies by \cite[Corollary~8.3.2]{GoLoan}:
\begin{equation}\tag*{\qed}
  \Error(\sigma_{\min}(M_q))\leq \|E\|
  <\ll n (\log N + D + n)\rr.
\end{equation}

\subsection{Bounding the error on the Newton step}

We next evaluate the error in the computation of
a Newton step. Our result is the following.

\begin{proposition}\label{coro:Newton}
There exists a universal constant $\bfe>0$ such that
given a system $q \in S(\Hd)$ and a point
$x \in S(\C^{n+1})$,
if the precision $u$ satisfies
$$
       u \leq \frac{\bfe}{ D^2\mun^2(q,x) n (D +\log N + n^2) },
$$
then the error $\Error (N_{q}(x))$ satisfies
$$
   \frac{\| \Error (N_{q}(x)) \|}{\| N_{q}(x) \|} \leq
    \frac{C (1 - \e) }{ 2 \pi (1+\e)^2 D^{3/2}\mun(q,x)},
$$
where $C$ and $\e$ are the constants introduced in
Proposition~\ref{prop:apps}.
\end{proposition}

We compute $N_q(x) - x$ by solving
the linear system
\begin{equation*}%\label{eq:ComputeNewton}
  \left[\begin{matrix} Dq(x)\\
  x^*  \end{matrix} \right] y = \left( \begin{matrix} q(x) \\
    0 \end{matrix} \right).
\end{equation*}
We denote $D_q =\left[\begin{matrix} Dq(x)\\
  x^*  \end{matrix} \right]$.

Recall the following result from~\cite[Chapter 7]{Higham96}
(in fact, Theorem 7.2 therein applied to $f = b / \|b\|$
and $E = A/ \|A\|$).

\begin{lemma}\label{lemma:LinearSystem}
Given a linear system $A x = b$, approximations $A'$ of $A$ and $b'$ of $b$
such that $\|A - A'\| \leq \epsilon$, $\|b - b'\| \leq \epsilon$ and $\epsilon \|A^{-1}\| < 1$, the
solution $x'$ of the perturbed system $A' x' = b'$ satisfies :
\begin{equation}\tag*{\qed}
\| x'-x \| \leq
\frac{\epsilon  \|A^{-1} \|}{1 - \epsilon  \|A^{-1} \|} (1 + \|x\|).
\end{equation}
\end{lemma}

Furthermore, from~\cite[Chapter 18]{Higham96}, the solution $\hat x$
of a linear system $Ax = b$,
where $A$ is non-singular, computed with a QR factorization
with finite precision satisfies
\begin{equation}\label{eq:linearSys}
(A + \Delta A) \hat x = b + \Delta b,
\end{equation}
where $|\Delta A| \leq n^2 \gamma_{cn} G |A|$, $|\Delta b| \leq n^2 \gamma_{cn} G |b|$, $\|G\|_F = 1$, $c \in \R$. Here, $|A|$ denotes the matrix with entries 
$|a_{ij}|$ and the same for $|\Delta A|$, $|b|$, and $|\Delta b|$.

\begin{lemma}\label{lemma:ErrorNewtonStep}
Let $q \in S(\mathcal{H}_d)$ and
$x\in S(\C^{n+1})$.
We can compute $N_q (x)$ with finite precision such that
$$
\frac{\|\Error (N_q(x)) \|}{\| N_q(x) \|} = \ll \mun(q,x) \cdot n \sqrt{D}
(D + \log N + n^2) \rr.
$$
\end{lemma}

\proof
To simplify notations, in this proof, we write $q(x)$ instead of $\left( \begin{matrix} q(x) \\
0 \end{matrix} \right)$.
Let $y$ denote our computed solution of $D_q y = q(x)$.

From~\eqref{eq:linearSys},
$\fl(y)$ is the solution of a system $(\Delta D_q + \fl(D_q)) y = \Delta q(x) + \fl(q(x))$  with $\|\Delta D_q \| = \ll n^3 \sqrt{D} \rr$
and $\|\Delta q(x) \| = \ll n^3\rr$.

From Proposition~\ref{prop:Errq},
given $q \in S(\mathcal{H}_d)$ and $x \in S(\C^{n+1})$, we
can compute $q(x)$ with finite precision $u$ such that $\|\Error
(q(x))\| = \ll D + \log N\rr$.
Obviously, $\|\Error(D_q)\|$ is not greater than the bound we 
computed for $\| \Error (Dq(x)H)\|$ in Proposition~\ref{prop:ErrD}.
Hence, $\Error(D_q) =\ll n \sqrt{D} (\log N + D + \log n)\rr$.

Furthermore, from~\eqref{eq:linearSys}, $\|\Delta D_q\| = \ll n^3 \sqrt{D} \rr$ and
$\|\Delta q(x)\| = \ll n^3\rr$.

Finally, both terms $\|\Delta D_q + \Error(D_q)\|$ and $\|\Delta q(x) + \Error( q(x))\|$ can be bounded by an expression
$\epsilon = \ll n \sqrt{D} (\log N + D + n^2)\rr.$

Thus, from Lemma~\ref{lemma:LinearSystem}, the error on $y$ is bounded as
\begin{eqnarray*}
\Error(y) & \leq &
       \frac{\epsilon  \|{D_q}^{-1} \|}
       {1 - \epsilon  \|{D_q}^{-1} \|}
       (1 + \|y\|)\\
       & = &  (1 + \|y\|) \ll \mun(q,x)  n \sqrt{D}
(D + \log N + n^2) \rr,
\end{eqnarray*}
the last line since $\|D_q^{-1}\| = \Oh(\mun(q,x))$ and $\frac{1}{1-\gamma_k}\leq \gamma_{4k}$.

Since $N_q(x)$ belongs to the tangent space to the unit sphere
$T_x$, $\|N_q(x)\| \geq 1$, and
\begin{equation*}
  1+ \|y\| = 1+ \|N_q(x) - x\| \leq 1 +  \|N_q(x)\| + \| x\| \leq 3 \|N_q(x)\|.
\end{equation*}

Hence,
$$
    \|\Error(y) \| =
    3  \|N_q(x)\| \ll \mun(q,x)  n
    \sqrt{D} (D + \log N + n^2) \rr.
$$
Then, the computation of $N_q(x)$ from $y=N_q(x) - x$ is a
simple addition and does not change the order of the errors.
\eproof

The proof of Proposition~\ref{coro:Newton} is now immediate.

\subsection{Bounding the error for $\Delta \tau$}

We evaluate here the errors in the computation of the quantity
$\Delta \tau$, that is, the size of
the current step in the homotopy.

\begin{proposition}\label{coro:DeltaTau}
For $x \in S(\C^{n+1})$, and $f,g,q \in S(\Hd)$ such that
$\dsp(f,g) \leq \pi/2$ define the quantity
$$
\Delta \tau := \frac{\lambda}{\dsp(f,g) D^{3/2}\mun^2(q,x)}.
$$

There exists a universal constant $\bff>0$
such that
\begin{equation}\label{eq:precisionDeltaTau}
   u \leq \frac{\bff}{ n (\log N + D + n) \mun^2(q,x)}
\end{equation}
implies
$$
   \Error(\Delta \tau)  \leq \frac{1}{4} \Delta \tau.
$$
\end{proposition}

To prove this proposition we rely on the following lemma.

\begin{lemma}\label{lemma:ErrorSigma2}
Given $x \in S(\C^{n+1})$ and $q \in S(\Hd)$, one
can compute $\sigma_{\min}^2(M_q)$ with finite precision $u$
such that
$$
  \Error ( \sigma_{\min}^2(M_q)) = \ll n(\log N + D + n) \rr.
$$
\end{lemma}

\proof
By Proposition~\ref{prop:sigmaM}, $\Error (\sigma_{\min}(M_q))
= \ll n(\log N + D + n) \rr$. Hence, we have
\begin{eqnarray*}
  |\fl( \sigma_{\min}^2(M_q)) -  \sigma_{\min}^2(M_q)|
  &\leq& 2 |\sigma_{\min}(M_q)| \ll n(\log N + D + n) \rr
  + \ll n(\log N + D + n) \rr ^2 \\
 &\leq& \ll n(\log N + D + n) \rr + \ll n(\log N + D + n) \rr,
\end{eqnarray*}
since, by Proposition~\ref{prop:ErrM},
$|\sigma_{\min}(M_q)| \leq
\|M_q\| \leq 1$. Thus,
\begin{equation}\tag*{\qed}
     \Error (\sigma_{\min}^2(M_q)) = \ll n(\log N + D + n) \rr.
\end{equation}

\proofof{Proposition~\ref{coro:DeltaTau}}
One has
\begin{eqnarray*}
 \fl(\Delta \tau)
 &=& \fl\left(\frac{\lambda}{\alpha D^{3/2} \mun^2 (q,x)}\right)\\
 &=& \frac{\lambda}{D^{3/2}} \fl(\frac{\sigma_{\min}^2(M_q)}{\alpha })
  (1 + \theta_{\Oh(1)})\\
 &=& \frac{\lambda \fl\left(\sigma_{\min}^2(M_q)\right)}{\alpha D^{3/2}}
  (1 + \theta_{\Oh(\log N)}),
\end{eqnarray*}
the last equality being from Lemma~\ref{lemma:alpha}, and thus by Lemma~\ref{lemma:ErrorSigma2}
\begin{eqnarray*}
 \Error(\Delta \tau)
 &=& \frac{\lambda}{\alpha D^{3/2}} (\ll n(\log N + D + n) \rr + \ll \log N \rr)\\
 &=& \frac{\lambda}{\alpha D^{3/2}} \ll n(\log N + D + n) \rr.
\end{eqnarray*}
%The error in the addition of $\Delta \tau$ and $\tau$ does not change
%this order, and thus, the same expression holds for $\Error( \tau +
%\Delta \tau).$

If $u$ satisfies (\ref{eq:precisionDeltaTau}) with a value of $\bff$ small
enough, the term $\ll n(\log N + D + n) \rr$ may be bounded by
$$
     \ll n(\log N + D + n) \rr \leq \frac{1}{4 \mun^2 (q,x)},
$$
and consequently
\begin{equation}\tag*{\qed}
 \Error(\Delta \tau) \leq
 \frac{\lambda }{4 \alpha D^{3/2} \mun^2(q,x)}
 = \frac{1}{4} \Delta \tau.
\end{equation}

\subsection{Bounding the distance between $\tilde q_\tau$ and
$q_\tau$}~\label{section:distqTau}

We evaluate here the error in the computation of $q_\tau$,
given $f,g$, and $\tau$.

\begin{proposition}~\label{prop:ErrQtau}
There exists a universal constant $\bfg$ such that the following holds.
Let $f,g \in S(\Hd)$ with $\dsp(f,g) \leq \frac{\pi}{2}(1+1/6)$
be given with roundoff error $u$. Let $\tau \in [0,1]$.
Then for all $A \in (0,1)$,
$$
u \leq \frac{\bfg \cdot A}{\log N}
$$
implies
$$
\|\fl (q_\tau) - q_\tau \| \leq A.
$$
\end{proposition}

We first bound the distance between the points $t f + (1-t)g$ and $t \fl(f) + (1-t)\fl(g)$,
without taking into account the error in the computation of $t$.

\begin{proposition}\label{prop:qtildeq}
Assume that $f,g,\tilde f, \tilde g \in S(\Hd)$ are such that $\dsp( f, g)\leq
\frac{\pi}{2}(1+1/60)$ and $\|f-\tilde f\|\leq
1/60$, $\|g-\tilde g\|\leq 1/60$.
For $t \in [0,1]$ define $q= t f + (1-t) g$ and $\tilde q = t \tilde f + (1-t) \tilde g$.
Then
$$
   \dsp(q, \tilde q)\leq
   2 \max\left\{\|f-\tilde f\| ,\|g- \tilde g\|\right\} .
$$
\end{proposition}

To prove Proposition~\ref{prop:qtildeq} we rely on the following
lemmas.

\begin{lemma}\label{lem:OptDist}
Let $f,g\in \Hd$ with  $\|f\|,\|g\|\geq \a>0$, $\|f-g\|\leq \b$
with $\a\geq \b/2$. Then
\begin{equation}\label{eq:opt}
\frac{|\<f,g\>|}{\|f\|\|g\|}\geq 1-\frac{\b^2}{2\a^2}.
\end{equation}
\end{lemma}

\proof
Pick any $f,g\in \Hd$ with $\|f\|,\|g\|\geq \a$, $\|f-g\|\leq
\b$, denote $r = (f+g)/2$, and let $s\in \Hd$ be such that
$\|s\|=\|f-g\|/2$, $s\perp r$. Then from the orthogonality of $r$
and $s$ we have
$$
\|r+s\|^2 = \|r-s\|^2 = \|r\|^2+\|s\|^2 =
\frac{\|f\|^2+\|g\|^2}{2}\geq \|f\|\|g\|\geq \a^2;
$$
also,
$$
\|(r+s)-(r-s)\| = 2\|s\| = \|f-g\|\leq \b.
$$
Therefore,
$$
\frac{|\<r+s,r-s\>|}{\|r+s\|\|r-s\|}\le
\frac{|\|r\|^2-\|s\|^2|}{\|f\|\|g\|}=
\frac{|\|f+g\|^2-\|f-g\|^2|}{4\|f\|\|g\|}\leq \frac{|\<f,g\>|}{\|f\|\|g\|}.
$$
Since $\|r+s\|=\|r-s\|$, we have
\begin{equation}\label{eq:opt01}
\min_{\|f\|,\|g\|\geq \a \atop \|f-g\|\leq
\b}\frac{|\<f,g\>|}{\|f\|\|g\|}= \min_{\|f\|=\|g\|\geq \a \atop
\|f-g\|\leq \b}\frac{|\<f,g\>|}{\|f\|\|g\|}.
\end{equation}

Now assume that $f,g\in \Hd$ with $\|f\|=\|g\|\geq \a$,
$\|f-g\|\leq \b$. Let
\begin{align*}
f' =
\frac{\b}{2}\cdot\frac{f-g}{\|f-g\|}+\frac{\sqrt{4\a^2-\b^2}}{2}\cdot\frac{f+g}{\|f+g\|};\\
g' =
\frac{\b}{2}\cdot\frac{g-f}{\|f-g\|}+\frac{\sqrt{4\a^2-\b^2}}{2}\cdot\frac{f+g}{\|f+g\|}.
\end{align*}
It is not difficult to check that $\|f'\|=\|g'\|=\a$ and
$\|f'-g'\|=\b$. Moreover,
$$
\frac{|\<f',g'\>|}{\|f'\|\|g'\|} = 1 - \frac{\b^2}{2\a^2}\leq
\frac{\|f\|^2+\|g\|^2}{2\|f\|\|g\|}-
\frac{\|f-g\|^2}{2\|f\|\|g\|}\leq\frac{|\<f,g\>|}{\|f\|\|g\|} .
$$
Therefore,
\begin{equation}\label{eq:opt02}
\min_{\|f\|=\|g\|\geq \a \atop \|f-g\|\leq
\b}\frac{|\<f,g\>|}{\|f\|\|g\|}= \min_{\|f\|=\|g\|= \a \atop
\|f-g\|=\b}\frac{|\<f,g\>|}{\|f\|\|g\|}.
\end{equation}
From \eqref{eq:opt01} and \eqref{eq:opt02} we have
$$
\min_{\|f\|,\|g\|\geq \a \atop \|f-g\|\leq \b}\frac{|\<f,g\>|}{\|f\|\|g\|}
= \min_{\|f\|=\|g\|= \a \atop
\|f-g\|=\b}\frac{|\<f,g\>|}{\|f\|\|g\|}\geq 1 - \frac{\b^2}{2\a^2},
$$
which shows \eqref{eq:opt}.
\eproof
\medskip

\begin{lemma}\label{lem:prop_cos}
Let $f,g\in \Hd$ with
$\|f-g \|\leq \min\{\|f\|,\|g\|\}.$ Then
\begin{equation}\label{eq:prop_cos}
\dsp(f,g)< \frac{2}{\sqrt{3}}\cdot \frac{\|f-g
\|}{\min\{\|f\|,\|g\|\}}.
\end{equation}
\end{lemma}

From Lemma~\ref{lem:OptDist} we have
$$
\cos \dsp(f,g) = \frac{\<f,g\>}{\|f\|\|g\|}\geq
1-\frac{\b^2}{2\a^2},
$$
where $\b = \|f-g\|$, $\a = \min\{\|f\|, \|g\|\}$. From the Taylor
expansion for $\cos$ we obtain
$$
\cos \dsp(f,g) \leq 1 - \frac{ \dsp^2(f,g)}{2}+
\frac{\dsp^4(f,g)}{24},
$$
therefore
$$
\frac{\dsp^4(f,g)}{12}-\dsp^2(f,g)+\frac{\b^2}{\a^2}\geq 0.
$$
Solving the relevant quadratic equation for $\dsp^2(f,g)$, we have
\begin{equation}\label{eq:opt04}
\dsp^2(f,g) \in \left(-\infty,
6\left(1-\sqrt{1-\frac{\b^2}{3\a^2}}\right)\right]\cup
\left[6\left(1+\sqrt{1-\frac{\b^2}{3\a^2}}\right),\infty\right).
\end{equation}
By our assumption $\b/\a\leq 1$, therefore,
$$
6\left(1+\sqrt{1-\frac{\b^2}{3\a^2}}\right) >\pi^2,
$$
and the interval on the right-hand side of \eqref{eq:opt04} is
irrelevant (as $\dsp(f,g)\leq \pi$).  We have (using
$\b/\a\leq 1$ again)
$$
\dsp^2(f,g) \leq
6\left(1-\sqrt{1-\frac{\b^2}{3\a^2}}\right)=\frac{2}{1+\sqrt{1-\frac{\b^2}{3\a^2}}}
\cdot \frac{\b^2}{\a^2}<\frac{4\b^2}{3\a^2},
$$
which yields \eqref{eq:prop_cos}.\eproof

\begin{lemma}\label{lem:MaxDistSphere}
Let $f, g \in \Hd$, $\|g\|=\|f\| = 1$, and $\dsp(f,g)\leq \frac{\pi}{2}(1+\d)$. Then,
given $t \in [0,1]$, $q(t) = t f + (1-t) g$ satisfies
\begin{equation}\label{eq:norm_qtau}
\|q(t)\|  \geq \sqrt{1-\frac{\pi^2(1+\d)^2}{16}}.
\end{equation}
\end{lemma}
\proof Consider the function $\varphi:\R\to \R$ defined as
follows:
$$
\varphi (t) = \|q(t)\|^2 = \|g\|^2+2 t
\Re\<g,f-g\>+t^2\|f-g\|^2.
$$
Observe that $\min_{t\in \R} \varphi(t)$ is attained at
$$
t^* = -\frac{2\Re\<g,f-g\>}{2\|f-g\|^2}
=\frac{\|g-f\|^2+\|g\|^2-\|f\|^2}{2\|f-g\|^2}=\frac{1}{2},
$$
and $\varphi(t^*) = \frac{1}{4}\|f+g\|^2$. We then have
$$
\|q(t)\|^2  \geq \frac{1}{4}\|f+g\|^2 =1-\frac{\|f-g\|^2}{4}
\geq 1- \frac{\dsp^2(f,g)}{4} \geq 1- \frac{\pi^2(1+\d)^2}{16},
$$
which gives us \eqref{eq:norm_qtau}. \eproof

\proofof{Proposition~\ref{prop:qtildeq}}
Observe that
\begin{align*}
\|q - \tilde q\|  & = \|t f+(1-t)g-t \tilde f -
(1-t)\tilde g\|\\
& \leq t\|f-\tilde f\|+(1-t)\|g-\tilde g\|\\
& \leq \max\left\{\|f-\tilde f\|, \|g- \tilde g\|\right\} \leq \frac{1}{60}.
\end{align*}
From Lemma~\ref{lem:MaxDistSphere} applied with $\d=\frac{1}{6}$ we have
$$
\|q\|\geq \sqrt{1-\frac{(1+1/60)^2}{4}}> \frac{3}{5},
$$
and hence
$$
\|\tilde q\|\geq \|q\|- \|q - \tilde q\| >
\frac{3}{5}-1/60 = \frac{7}{12}.
$$
Now applying Lemma~\ref{lem:prop_cos}, we have
\begin{equation}\tag*{\qed}
\dsp(q, \tilde q) \leq \frac{2}{\sqrt{3}}\cdot
\frac{\|q - \tilde q\|}{\min\{\|q\|, \|\tilde
q\|\}} \leq \frac{2}{\sqrt{3}}\cdot \frac{
\max\left\{\|f-\tilde f\|, \|g- \tilde g\|\right\}
}{\frac{7}{12}} \leq 2 \max\left\{\|f-\tilde f\|, \|g- \tilde g\|\right\}.
\end{equation}
%Observe that when $\d<\frac{1}{6}$, the denominator is greater
%than $\sqrt{3}$, hence, we have
%\begin{equation}\tag*{\qed}
%\dsp(q_\tau, \tilde q_\tau) \leq \frac{8}{3}\max\left\{\|f-\tilde
%f\| ,\|g- \tilde g\|\right\}.
%\end{equation}

\proofof{Proposition~\ref{prop:ErrQtau}}
Let us denote $\tilde f = \fl(f)$, $\tilde g = \fl(g)$ and $\tilde t = \fl (t)$.
Let $\tilde q_\tau $ denote $\fl(q_\tau)$ and $\hat q_\tau$ the system $t \tilde f + (1-t) \tilde g.$
%Since each coefficient of $f$ is computed with
%roundoff error $u$,
%$$
%\| \tilde f - f \| =\left( \sum_\a (f_\a (1+ \t_1) - f_\a)^2 \right)^{1/2} = \t_{1} \left(\sum_\a {f^2}_\a\right)^{1/2} = \t_{1}\leq 2u,
%$$
By hypothesis, both $\| f-\tilde f\|$ and $\|g - \tilde g\|$ are not greater than $u$.

Thus, by Proposition~\ref{prop:qtildeq}, if $u \leq 1/60$,
$$
\|\hat q_\tau - q_\tau \| \leq 2 u.
$$

From Proposition~\ref{prop:Errt}, $\tilde t  = t (1 + \theta_{\Oh(\log N)}).$
Thus, there exists a constant $\bfg$ such that for all $A\in (0,1)$, $u \leq \frac{\bfg A}{\log N}$ implies
$$
\|\hat q_\tau - \tilde q_\tau \| \leq A/2.
$$

Taking $\bfg \leq 1/60$, $u \leq \frac{\bfg A}{\log N}$ ensures

\begin{equation}\tag*{\qed}
\|\tilde q_\tau - q_\tau\| \leq \|\tilde q_\tau - \hat q_\tau \| + \|\hat q_\tau - q_\tau \| \leq \frac56 A < A.
\end{equation}

\subsection{Estimates for $u$}\label{section:initU}

Along our homotopy, the precision needed to guarantee correctness 
varies with the system $q$ considered. In our variable-precision algorithm 
we will want to keep this precision at all times
within the interval $[ \frac12 \bfB(q,x), \bfB(q,x)]$
with
\begin{equation}\label{eq:bfB}
   \bfB(q,x) := \frac{k_2}{n D^{2} (\log N + D + n^2)
   \mun^2(q,x)},
\end{equation}
where $k_2$ is a universal positive constant (that will be specified in
Definition~\ref{def:k2}).

Now, since $q$ and $x$ vary at each iteration, one has to
update the precision as well. To do so one faces an obstacle.
When computing $u$
%---for instance, as $u:=\frac34\bfB(q,x)$ in the
%{\tt repeat} loop---
we actually obtain a quantity
$\fl(u)$ which depends on the current precision
and this current precision has been computed
in the previous iteration. Proposition~\ref{prop:repeat_u}
below shows that this obstacle can be overcome.

\begin{proposition}\label{prop:repeat_u}
If $u\leq (1+ \e)^6 \bfB(q,x)$ then
$$
\Error \left(\frac{3}{4} \bfB(q, x)\right)
\leq \frac{1}{4} \bfB(q, x).
$$
In particular, when computing $u:=\frac34\bfB(q,x)$
the computed quantity satisfies
$\fl(u)\in [\frac12\bfB(q,x),\bfB(q,x)]$.
\end{proposition}

Towards the proof of the proposition above we define
$$
    B(q,x) := \frac{1}{n D^{2} (\log N + D + n^2) \mun^2(q,x)}
$$
so that $\bfB(q,x)=k_2B(q,x)$. Our first lemma bounds
the error in the computation of $B(q,x)$.

\begin{lemma}\label{lemma:computeBq}
There exists a positive universal constant $k_3$ such that the following is true.
Assume $u \leq \frac{1}{2 k_3 n(\log N + D + n)}$. Then
for any $x \in S(\C^{n+1})$  and $q \in S(\Hd)$,
one can compute $B(q,x)$ such that
$$
   \Error (B(q,x)) \leq \frac{k_3 u }{D^2}. 
$$
\end{lemma}

\proof
From Lemma~\ref{lemma:ErrorSigma2}, we can compute
$\sigma_{\min}^2(M_q)$ with error $\ll n (\log N + D + n) \rr$; we
can compute $B(q,x)$ such that
$$
    \fl(B(q,x)) = \frac{\fl (\sigma_{\min}^2(M_q))}{n D^{2}
    (\log N + D + n^2)} (1 + \t_6),
$$
and thus
$$
    \Error (B(q,x)) = \frac{\ll n (\log N + D + n^2) \rr }
    { n D^2 (\log N + D + n^2)}.
$$
It follows that there exists a constant $k_3$ such that
$$
 \Error (B(q,x)) = \frac{\theta_{k_3 n (\log N + D + n^2)}}
{ n D^2 (\log N + D + n^2)}.
$$

Thus, when $u \leq \frac{1}{2k_3 n(\log N + D + n^2)}$, the denominator in
$\gamma_{k_3 n (\log N + D + n)}$ is greater than $1/2$, and
\begin{equation}\tag*{\qed}
\Error (B(q,x)) \leq \frac{k_3 u}{D^2}.
\end{equation}

\begin{corollary} \label{coro:IterateU}
Let $k_2$ be a positive constant such that
$k_2\leq\frac{1}{4 k_3 (1+ \e)^6}$.
The condition $u \leq (1 + \e)^6 k_2 B(q,x)$ ensures that
$$
  \Error \left(\frac{3}{4} k_2 B(q, x)\right) \leq \frac{1}{4} k_2 B(q, x).
$$
\end{corollary}

\proof
If $u$ is less than or equal to $(1 + \e)^6 B(q,x)$, since $\mun (q, x)$ is always greater than 1, if we choose $k_2$ not greater than $\frac{1}{2 k_3(1+ \e)^6}$, $u$ will be less than or equal to $\frac{1}{2 k_3 n(\log N + D + n^2)}$.
Thus, one has
$$
\Error (B(q,x)) \leq \frac{k_3 u}{D^2}
\leq  \frac{k_2 k_3 (1 + \e)^6}{D^2} B(q, x).
$$
Taking $k_2 \leq \frac{1}{4 k_3 (1+ \e)^6}$ one has
$\Error (\frac{3}{4} k_2 B(q, x))\leq \Error (k_2 B(q, x))
\leq \frac{1}{4} k_2 B(q, x)$.
\eproof

We now have all the conditions that the constant $k_2$ must fulfill.

\begin{definition}~\label{def:k2}
Let $k_2$ be a positive constant chosen small enough such that
%$u \leq B(q,x)$ ensures, that

\begin{enumerate}
\item[(i)] $k_2 B(q,x)$ is smaller than the bound on $u$ in
Proposition~\ref{coro:DeltaTau},

\item[(ii)] $k_2 B(q,x) (1+\e)^4$ is smaller than the bound on $u$ in
Proposition~\ref{coro:Newton},

\item[(iii)] $k_2 \leq \frac{\bfg C (1-\e)}{12(1+\e)^5}$,
where $\bfg$ is defined in Proposition~\ref{prop:ErrQtau}.

\item[(iv)] $k_2$ verifies the condition in Corollary~\ref{coro:IterateU}.

\end{enumerate}
\end{definition}

The first three conditions ensure that the precision will be good
enough for the computation of the values of $\Delta \tau$, of
the Newton operator and of $q_\tau$.
The fourth condition is needed for the computation of
$B(q,x)$ itself.
\medskip

\proofof{Proposition~\ref{prop:repeat_u}}
From~\eqref{eq:bfB}, the bound
$\bfB(q,x)$ equals $k_2 B(q,x)$
and the result now follows from
Corollary~\ref{coro:IterateU}.
\eproof

%\begin{corollary}
%At each step in the {\tt repeat} loop of \ALHF\ the computed value
%for $u$ is in $[\frac12 \bfB(q_\tau,x),\bfB(q_\tau,x)]$.
%\end{corollary}

%\proof
%ADD IT.
%\eproof

\section{Analysis of the Homotopy}\label{sec:homotopy}

We next describe with more detail our procedure \ALHF\
---{\bf A}daptive {\bf L}inear {\bf H}omotopy
with {\bf F}inite precision--- to follow the path
$\{(q_\tau,\z_\tau)\mid \tau\in[0,1]\}$.

All the certifications on an execution of \ALHF\ will be
for inputs satisfying certain conditions. We thus define the
notion of admissible input for \ALHF.

\begin{definition}\label{def:admInput}
An \emph{admissible input} for algorithm \ALHF\ consists
of
\begin{itemize}
\item A function {\tt read\_input$_f$(\ )}, that returns an approximation
of a system $f \in S(\Hd)$ with the current round-off unit. That is,
the instruction {\tt read\_input$_f$(\ )} returns a system $f'$
such that the coefficients $a'_\alpha$ of the polynomials $f'_i$ satisfy
$$
| a'_\alpha - a_\alpha | \leq u |a_\alpha|,
$$
where $a_\alpha$ is the coefficient of the monomial of the
same degree $\alpha$ of $f_i$. In particular, this implies that
$$
\| f - f' \| \leq u \|f \|.
$$

Note that {\tt read\_input$_f$(\ )} is not required to be computable.

\item An auxiliary system $g \in S(\Hd)$, supposed to be given
exactly.

\item An approximate zero $x\in S(\C^{n+1})$ of $g$ satisfying
$$
\dpr(\z,x) \leq\frac{C}{D^{3/2}\mun(g,\zeta)}
$$
for its associated zero $\zeta$.

\item An initial round-off unit $u \in \R_+$ such that
$$
  u \leq \bfB(g,x).
$$
\end{itemize}

For clarity, we denote such a tuple $(f,g,x,u)$ and we
refer to it as an input to \ALHF\ even though $f$ is not
given directly and the precision $u$ is not passed as a
parameter (it is a global variable in \MDF).
\end{definition}

Define $\lambda:=\frac{2C(1-\e)}{5(1+\e)^4} \approx 5.37 \cdot 10^{-3}$.
\medskip

\begin{center}
\algo
\> Algorithm \ALHF\\[2pt]
\>{\bf input} $(f,g,x)$\\[2pt]
\>\>$\tilde f := {\tt read\_input}_f(\ )$\\[2pt]
%\>\>$\tilde f:=\frac{\tilde f}{\|\tilde f\|}$\\[2pt]
\>\> if $\dsp(\tilde f,g)\geq \frac{\pi}{2}$ then
 $g:=-g$\\[2pt]
\>\>$\tau:=0$, $\tilde q_\tau:=g$\\[2pt]
\>\> {\tt repeat}\\[2pt]
\>\>\>$\Delta\tau:= \frac{\lambda}{\dsp(\tilde f,g) D^{3/2}\mun^2(\tilde q_\tau,x)}$\\[2pt]
\>\>\>$\tau:=\min\{1,\tau+\Delta\tau\}$\\[2pt]
%\>\>\>$\eta:=\frac{C(1-\e)}{12 (1+\e)^5 D^{3/2}
%         \mun^2(\tilde q_\tau,x)}$\\[2pt]
\>\>\>$\tilde f := {\tt read\_input}_f(\ )$\\[2pt]
\>\>\>$\tilde q_\tau:=t(\tau)\tilde f+(1-t(\tau)) g$\\[2pt]
\>\>\>$\tilde q_\tau:=\frac{\tilde q_\tau}{\|\tilde q_\tau\|}$\\[2pt]
\>\>\>$x:=N_{\tilde q_\tau}(x)$\\[2pt]
\>\>\>$x: = \frac{x}{\|x\|}$\\[2pt]
\>\>\>$u:= \frac34 \bfB(\tilde q_\tau,x)$\\[2pt]
\>\> {\tt until} $\tau= 1$\\[2pt]
\>\> {\tt RETURN $x$}
\falgo
\end{center}
\medskip

\begin{remark}
The algorithm \ALHF\ is a finite-precision adaptation
of the algorithm \ALH\ in~\cite{BC09}. It has a
slightly smaller stepsize parameter
$\lambda$. By the parameter $f$ given to \ALHF, we mean,
that the algorithm is given as input the procedure
{\tt read\_input$_f$} that returns finite precision approximations of $f$.
\end{remark}

We may use \ALHF\ to define a finite precision version \MDF\ of \MD.

\algo
\> Algorithm \MDF\\[2pt]
\>{\bf input} $f\in \Hd$ \\[2pt]
\>\> $u:=\frac{k_2}{n D^2( \log N + D +n^2) 2 (n+1)^D}$\\[2pt]
\>\> {\tt run \ALHF\ on input $(f,\oU,\bz_1)$} \falgo
\medskip

%A main step to bound the average expected cost of \LVF\
%is to estimate the number of operations
%$\cost_{\ALHF}(f,\bar{w})$ performed by \ALHF\
%on a given execution. Here $f$ denotes the input and
%$\bar{w}$ the sequence of points in $W$ drawn during
%this execution.  A similar remark applies to the
%average expected finest precision $u_*$ used by \LVF.
%The following theorem provides bounds for in
%this main step. The following function on
%homotopy paths will be used.

To a pair $f\in S(\Hd)$ and $(g,\zeta)\in V_{\P}$
we associate the number
$$
   \mu_*(f,g,\z):=
   \max_{\tau\in[0,1]}\mun(q_\tau,\z_\tau).
$$

%% IB à reformuler
\begin{theorem}\label{thm:main1}
Let $(f,g,x,u)$ be an admissible input of \ALHF. Then:
\begin{enumerate}

\item[(i)]
If the algorithm \ALHF\ stops on input $(f,g,x)$,
it returns an approximate zero of $f$.

\item[(ii)]
Assume \ALHF\ stops on input $(f,g,x)$. Then,
the number of iterations
$K(f,g,x)$ performed by \ALHF\ satisfies
$$
  K(f,g,x) \leq B(f,g,\z)+B(-f,g,\z)
$$
where
$$
 B(f,g,\z):=408\, \dsp(f,g) D^{3/2}\int_0^1
  \mun^2 (q_\tau,\z_\tau)d\tau.
$$
Consequently the number of performed arithmetic operations
$\cost_{\ALHF}(f,g,x)$ is bounded by
$$
  \cost_{\ALHF}(f,g,x)\leq \Oh(N)\big(B(f,g,\z)+B(-f,g,\z)\big).
$$
If \ALHF\ does not stop then either $B(f,g,\z)$ or $B(-f,g,\z)$ is unbounded,
and either the segment $E_{g,f}$ or $E_{g,-f}$ intersects $\Sigma$.

\item[(iii)] Furthermore, the finest precision $u_*(f,g,x)$ required
during the execution is bounded
from below by
$$
  u_*(f,g,x)=\Omega \left( \frac{1}{n D^{2}
    (\log N + D + n^2)\mu_*^2(f,g,\z)}\right).
$$
\end{enumerate}
\end{theorem}

\subsection{Bounding errors in the homotopy}

We begin with a simple consequence of
Proposition~\ref{prop:apps}.

\begin{proposition}\label{prop:cdx}
Assume $D \geq 2$. Let $p_0,p_1\in S(\Hd)$, let $\zeta$ be a zero of
$p_0$, and $A$ a positive constant not greater than $C$
such that
$$
\dsp (p_0,p_1) \leq \frac{A}{(1+\e) D^{3/2} \mun^2(p_0,\z)}.
$$
Then the path $E_{p_0,p_1}$ can be lifted to a path in $\Vp$
starting in $(p_0,\z)$. In addition, the zero $\chi$ of $p_1$
in this lifting satisfies
$$
  \dpr (\z, \chi) \leq \frac{A}{D^{3/2} \mun(p_1,\chi)}.
$$
Finally, for all $p_\tau\in E_{p_0,p_1}$, if $\z_\tau$ denotes
the zero of $p_\tau$ in this lifting, we have
$$
  \frac{1}{1+\e}\mun(p_0,\z)\leq \mun(p_\tau,\z_\tau)
  \leq (1+\e)\mun(p_0,\z).
$$
\end{proposition}

\proof
For each $\tau \in [0,1]$, let $p_\tau$ be the point of the
segment $[p_0,p_1]$ such that
$\dsp (p_0, p_\tau) = \tau \dsp (p_0,p_1)$.

Let $\tau_*$ be such that
$\int_{0}^{\tau_*}\mun(p_\tau,\zeta_\tau)\|\dot{p}_\tau\|d\tau
=\frac{A}{D^{3/2}\mun(p_0,\z)}$, or
$\tau_*=1$, or the path $E_{p_0,p_1}$ cannot be lifted to $V$
beyond $\tau_*$, whichever is the smallest.
Then, for all $\tau \in[0,\tau_*]$, using that
$\|\dot{\z}_\tau\| \le \mun(p_\tau,\z_\tau)\, \|\dot{p}_\tau\|$
(cf.~\cite[\S12.3-12.4]{bcss:95}) we have
\begin{eqnarray*}
 \dpr(\zeta,\zeta_\tau) &\leq &
    \int_{0}^\tau\|\dot{\z}_s\|\, ds
    \leq \int_{0}^{\tau_*}\mun(p_s,\z_s)\,
    \|\dot{p}_s\| ds \\
    &\leq& \frac{A}{D^{3/2}\mun(p_0,\z)}.
\end{eqnarray*}
It is therefore enough to show that $\tau_* = 1$.
Suppose to the contrary, that $\tau_* < 1$.

Since $\mun(p_\tau, \zeta_\tau) \geq 1 $, for every $\tau$,
\begin{equation*}
 \dsp(p_0,p_\tau) \leq\dsp(p_0,p_1)\leq
 \frac{A}{D^{3/2}\mun(p_0,\z)}.
\end{equation*}

Since $A\leq C$ the bounds on $\dsp(p_0,p_\tau)$ and
$\dpr(\z,\z_\tau)$ allow us to apply Proposition~\ref{prop:apps}
and to deduce, for all $\tau\in[0,\tau_*]$,
\begin{equation}\label{eq:cdx1}
    \frac{\mun(p_0,\z)}{1+\e} \leq
    \mun(p_\tau,\zeta_\tau)
    \leq (1+\e) \mun(p_0,\z).
\end{equation}

We have
\begin{eqnarray*}
  \frac{A}{D^{3/2}\mun(p_0,\z)}
  & = &  \int_{0}^{\tau_*}
    \mun(p_\tau,\z_\tau)\|\dot{p}_\tau\| d\tau
    \quad \mbox{(by definition of $\tau_*$)}\\
   & \leq\ &   (1+\e)\mun(p_0,\z) \int_{0}^{\tau_*}\|\dot{p}_\tau\| d\tau \quad \mbox{(by (\ref{eq:cdx1}))}\\
   &=&\  \dsp(p_0,p_{\tau_*})(1+\e)\mun(p_0,\z),
\end{eqnarray*}
and thus
$$
\dsp(p_0,p_{\tau_*}) \geq \frac{A}{(1+\e) D^{3/2}\mun^2(p_0,\z)} \geq \dsp(p_0,p_1),
$$
which leads to a contradiction with $\tau_* < 1$, and finishes the proof.
\eproof

The next proposition puts together many of the results obtained
thus far. The general idea for its proof
closely follows~\cite[Theorem~3.1]{BC09}
(which in turn is a constructive version of the
main result in~\cite{Bez6}) making
some room for errors.
\medskip

%Fix $f,g\in\Hd$, $x\in\C^{n+1}$, and $u\in(0,1)$.
%% and $\z\in\C^{n+1}$ such that $(g,\z)\in V$. Assume
%%that the path $E_{g,f}$ does not meet the discriminant
%%variety and hence, it can be lifted to a path $(q_\tau,\z_\tau)$ on $V$.
%%
%Suppose, that there exists a zero $\zeta$ of $g$ such that
%$\dpr(\z,x) \leq\frac{C}{D^{3/2}\mun(g,\zeta)},$ and
%suppose also $ u \leq \bfB(g,x)$.
Let $(f,g,x,u)$ be an admissible input for algorithm \ALHF.

Let $0=\overline{\tau_0}<\overline{\tau_1}
<\overline{\tau_2}<\ldots$,
$x=\overline{x_0},\overline{x_1},\overline{x_2},\ldots$, and
$\overline{u}_0,\overline{u}_1,\overline{u}_2,\ldots,$
be the sequences of $\tau$-values, points in $S(\C^{n+1})$
and precisions generated by the algorithm \ALHF\ on
the admissible input $(f,g,x,u)$.
Let $\tilde f_i$ be the approximation of the input $f$ on the $i$th iteration.

Let $E_{g,f}$ be the path
with endpoints $g$ and $f$. To simplify notation we write
$q_i$ instead of $q_{\overline{\tau_i}}$
and $\z_i$ instead of $\zeta_{\overline{\tau_i}}$.
Similarly, we denote by $\tilde{q}_i$ the computed approximation
of $q_i$ ---that is,
$\tilde q_i=\fl\big(t(\tau_i)\tilde f_i +(1-t(\tau_i)) g\big)$---,
by $x_{i+1}$ the exact value of
$N_{q_i}(\overline{x_i})$, and
by $\tau_{i+1}$ the exact value of
$\overline{\tau_i} + \Delta \tau$.
%, and by $\overline{x_{i+1}}$ and $\overline{\tau_{i+1}}$ the
%values for $x_{i+1}$ and $\tau_{i+1}$
%actually computed by \ALHF.

\begin{proposition}~\label{prop:Recurence}
Let $(f,g,x,u)$ be an admissible input for \ALHF.
Let $k$ be the number of iterations of \ALHF\ on input $(f,g,x,u)$ ---
that is, either $k=\infty$ or $\tau_k=1$, $q_k = f$.
With the notations above, for all $i\in \{0, \ldots, k-1\}$,
the following inequalities are true:
\medskip

\noindent {\bf (a)}\quad $\displaystyle
 \dpr(\overline{x_i},\zeta_i)\leq \frac{C}{D^{3/2}\mun(q_i,\zeta_i)}$
\smallskip

\noindent {\bf (u)}\quad $\displaystyle
\frac{\bfB(\tilde q_i,\overline{x_i})}{2} \leq
\overline{u}_i \leq \bfB(\tilde q_i,\overline{x_i})$
\smallskip

\noindent {\bf (x)}\quad $\displaystyle
  \dsp(q_i,\tilde{q}_i)\leq \frac{C(1-\e)}{12 (1+\e) D^{3/2}
\mun(q_i,\zeta_i)}$
\smallskip

%%   \noindent {\bf (y)}\quad $\displaystyle
%%    \forall \tau \in [0,1] \, \dpr(\zeta_\tau,\tilde{\zeta_\tau})\leq \frac{C(1-\e)}{ 12(1+\e)
%%    D^{3/2}\mun(q_i,\zeta_i)}$
%%   \smallskip

%
% \noindent {\bf (b)}\quad $\displaystyle
%     \frac{\mun(q_i,\overline{x_i})}{(1+\e)}\leq
%     \mun(q_i,\zeta_i)\leq (1+\e)\mun(q_i,\overline{x_i})$
% \smallskip

\noindent {\bf (c)}\quad $\displaystyle
 \dsp(q_i,q_{i+1})\leq \frac{(1-\e) C}{2 (1+\e) D^{3/2}
\mun(q_i,\zeta_i)}$
\smallskip

 \noindent {\bf (d)}\quad $q_{i+1}$ has a zero $\z_{i+1}$ such that
 $\displaystyle
  \dpr(\zeta_i,\zeta_{i+1})\leq \frac{(1-\e)C}{2 (1+\e)
  D^{3/2}\mun(q_i,\zeta_i)}$
 % \,\frac{(1-\e)}{(1+\e)}$
 \smallskip

\noindent {\bf (e)}\quad $\tilde q_{i+1}$ has a zero $\tilde \z_{i+1}$
such that
 $\displaystyle
 \dpr(\overline{x_i},\tilde{\zeta}_{i+1})\leq \frac{C \left( (1+ \e) + 7/12 (1- \e)\right)}{D^{3/2}
 \mun(q_{i+1},\zeta_{i+1})}$
\smallskip

\noindent
Inequalities {\bf (a)}, {\bf (u)}, and {\bf (x)} hold for $k$ as well
in case $k<\infty$.
\end{proposition}

Proposition~\ref{prop:Recurence} puts together all the needed
bounds to ensure the proper work of \ALHF. Statement ({\bf a},$i$)
ensures that $\overline{x_i}$ is ``close enough'' to $\zeta_i$. That
is, $\overline{x_i}$ is not just an approximate zero of $q_i$, but
also an approximate zero for polynomials in a certain neighborhood of
$q_i$ on $E_{g,f}$. Statements ({\bf c},$i$) and ({\bf d},$i$) show
that (taking into account computational errors) our step along the
homotopy is so small that the next polynomial $q_{i+1}$ belongs to
this neighborhood. We hence arrive at ({\bf e},$i$), which essentially
means that $\overline{x_i}$ is an approximate zero of $q_{i+1}$
associated with $\zeta_{i+1}$. Therefore, the Newton step (with
computational errors accounted for) brings the next iterate
$\overline{x_{i+1}}$ close enough to $\zeta_{i+1}$ to ensure that
({\bf a},$i+1$) holds again. Making sure that ({\bf u}) holds on every
iteration, we guarantee that computational errors are small enough to
allow all the other steps of the proof (({\bf a}), ({\bf c}), ({\bf
d}) and ({\bf e})) to be carried through.
\medskip

\proofof{Proposition~\ref{prop:Recurence}}
We proceed by induction by showing, that ({\bf a},$i$), ({\bf u},$i$)
and ({\bf x},$i$) imply successively ({\bf c},$i$), ({\bf d},$i$),
({\bf x},$i+1$), ({\bf e},$i$), and finally ({\bf a},$i+1$) and
({\bf u},$i+1$). %Note that the relation $w_i\prec w_{i+1}$
%implies that $\Gamma(w_i)\supset \Gamma(w_{i+1})$ and
%hence that induction hypothesis always hold for the set of
%considered pairs.

Inequalities ({\bf a}) and ({\bf u}), for $i=0$ hold
by hypothesis, and (({\bf x}), 0) is obvious since
$\tilde q_0 = q_0 = g$.

This gives us the induction base.
Assume now that ({\bf a}), ({\bf u}) and ({\bf x}) hold for some
$i \le k-1$.

We now show ({\bf c},$i$) and ({\bf d},$i$).

Observe that together with ({\bf a},$i$) and ({\bf x},$i$),
Proposition~\ref{prop:apps} implies
\begin{equation}\label{eq:b}
   \frac{\mun(\tilde{q_i},\overline{x_i})}{(1+\e)}\leq
   \mun(q_i,\zeta_i)\leq (1+\e)\mun(\tilde{q_i},\overline{x_i}).
\end{equation}

By ({\bf u},$i$) and Definition~\ref{def:k2} our
precision $\overline{u}_i$ satisfies~\eqref{eq:precisionDeltaTau}
for the pair $(\tilde{q_i},\bar{x_i})$.
Therefore, by
Proposition~\ref{coro:DeltaTau} and the definition of
$\Delta \tau$ in \ALHF\  we have
\begin{eqnarray*}
  \a(\overline{\tau_{i+1}}-\overline{\tau_i})
  &\leq& \a (\Error(\Delta \tau) + \tau_{i+1} -\overline{\tau_i})\\
  &\leq&\a\frac{5}{4}\Delta\tau
  \leq \frac{\lambda (1+\frac{1}{4})}{D^{3/2}
  \mun^2(\tilde{q}_i,\overline{x_i})}.
\end{eqnarray*}

So, using \eqref{eq:b} and since $\lambda:=\frac{2C(1-\e)}{5(1+\e)^4}$, we obtain
$$
  \dsp(q_i,q_{i+1}) \,= \,\a(\bar{\tau_{i+1}}-\bar{\tau_i})
  \leq \frac{C (1-\e)}{2(1+\e)^4 D^{3/2}
  \mun^2(\tilde{q}_i,\overline{x_i})}
  \leq \frac{C (1-\e)}{2(1+\e)^2 D^{3/2}
  \mun^2(q_i,\z_i)}.
$$
Since $\mun(q_i,\z_i)$ is always greater than or equal to 1, ({\bf c},$i$) holds, and ({\bf d},$i$)
is the direct consequence of Proposition~\ref{prop:cdx} applied to
$(q_i,q_{i+1})$ and $\zeta_i$, with $A=\frac{C (1-\e)}{2(1+\e)}$.

This application of Proposition~\ref{prop:cdx} furthermore
ensures that, for all $\tau\in[\bar{\tau_i},\bar{\tau_{i+1}}]$,
\begin{equation}\label{eq:btau}
    \frac{\mun(q_i,\zeta_i)}{1+\e} \leq
    \mun(q_\tau,\zeta_\tau)
    \leq (1+\e) \mun(q_i,\zeta_i),
\end{equation}
and, in particular,
\begin{equation}\label{eq:b1}
    \frac{\mun(q_i,\zeta_i)}{1+\e} \leq
    \mun(q_{i+1},\zeta_{i+1})
    \leq (1+\e) \mun(q_i,\zeta_i).
\end{equation}

%Since the current precision $\overline{u}_i$ satisfies
%(\ref{eq:precisionDeltaTau}) for the pair $(\tilde{q_i},\bar{x_i})$,
%by Proposition ~\ref{coro:DeltaTau}, the computed value of
%$\eta$ for this pair verifies
%$\Error (\eta ) \leq \eta/4$, and thus $\fl (3/4 \eta) \leq \eta$.
%It follows that ${\tt refine} (w_i,9\eta/32 )$ returns a
%triple $(w_{i+1},\tilde g_{i+1},\bar{x_{i+1}})$
%satisfying, for all $(g,\z)\in\Gamma(w_{i+1})$,
%$\|g-\tilde g_{i+1}\|\leq \frac{3}{8}\eta$. Similarly,
%since $\bar{u}_i\le \bfB(\tilde{q}_i,\bar{x_i}) \leq ??$
%we deduce from Proposition~\ref{prop:read_input} that
%${\tt read\_input} (\ )$ produces a $\tilde f$
%satisfying $\|f-\tilde f\|\leq \frac{3}{8}\eta$.

Since $u \leq \bfB(\tilde q_{i},\bar x_{i})$ and from Definition~\ref{def:k2},
we can apply Proposition~\ref{prop:ErrQtau}
with $A = \frac{C (1-\e)}{12 (1+\e)^5 D^{3/2} \mun^2(\tilde q_{i},\bar x_{i})}$
and we get
%We have that $\|\tilde f - f \| \leq u$. Thus,
%we may use Proposition~\ref{prop:qtildeq} with $\tilde g = g$
%to deduce
\begin{eqnarray}\label{eq:eta1}
  \dsp(q_{i+1},\tilde{q}_{i+1})
  &\leq& \frac{C (1-\e)}{12 (1+\e)^5 D^{3/2} \mun^2(\tilde q_{i},\bar x_{i})} \notag\\
 % \frac{8}{3}\max\{ \|f-\tilde f\|,0\}
 % \leq \frac{8}{3} \bfB (\tilde q_{i},\bar x_{i})
  % \eta = \frac{\frac{C(1-\e)}{12 (1+\e)^5}}{ D^{3/2} \mun^2(\tilde q_{i},\bar x_{i})}
  %\notag\\
  &\leq & \frac{\frac{C(1-\e)}{12 (1+\e)^3}}{D^{3/2}
    \mun^2(q_{i},\zeta_{i})} \quad \mbox{(from (\ref{eq:b}))}
\end{eqnarray}
and, hence, using (\ref{eq:b1}),
\begin{equation}\label{eq:x}
  \dsp(q_{i+1},\tilde{q}_{i+1})
  \leq  \frac{\frac{C(1-\e)}{12 (1+\e)}}{D^{3/2}
     \mun^2(q_{i+1},\zeta_{i+1})}.
\end{equation}
Since $\mun(q_{i+1},\z_{i+1})\geq 1$ this shows
({\bf x}, $i+1$).

We can now use ({\bf x}, $i$), ({\bf c}, $i$), and \eqref{eq:eta1}
to bound $\dsp(\tilde q_{i},\tilde{q}_{i+1})$ as follows,
\begin{eqnarray}\label{eq:bound_tildes}
\dsp(\tilde q_{i},\tilde{q}_{i+1})
&\leq& \dsp(\tilde q_{i},q_{i})
  + \dsp(q_{i},q_{i+1})
  + \dsp(q_{i+1},\tilde{q}_{i+1})\notag \\
&\leq& \frac{\frac{C(1-\e)}{12 (1+\e)}}{D^{3/2}
\mun(q_i,\zeta_i)}
+\frac{\frac{C(1-\e)}{2 (1+\e)}}{ D^{3/2}
\mun(q_i,\zeta_i)}
+\frac{\frac{C(1-\e)}{12 (1+\e)^3}}{D^{3/2}
    \mun^2(q_{i},\zeta_{i})} \notag\\
&<& \frac{\frac{C(1-\e)}{(1+\e)}}{D^{3/2}
\mun(q_{i},\zeta_{i})}\notag\\
&\leq& \frac{C(1-\e)}{D^{3/2}
\mun(\tilde q_{i},\bar x_{i})}
\end{eqnarray}
the third inequality using $\mun(q_{i},\z_{i})\geq 1$
and the last from~\eqref{eq:b}.
We can similarly bound distances between zeros and their
approximations. Indeed, using \eqref{eq:x},
Proposition~\ref{prop:cdx} applied
to $(q_{i+1}, \tilde{q}_{i+1})$ and $\z_{i+1}$,
with $A = \frac{C(1-\e)}{12}$,
ensures the existence of a zero $\tilde \z_{i+1}$
of $\tilde q_{i+1}$ such that
\begin{equation}\label{eq:y}
  \dpr(\zeta_{i+1},\tilde{\zeta}_{i+1}) \leq
   \frac{C(1-\e)}{ 12D^{3/2}\mun(q_{i+1},\zeta_{i+1})}.
\end{equation}
Next we use the triangle inequality to obtain
\begin{eqnarray*}
  \dpr(\overline{x_i},\tilde{\zeta}_{i+1})
  &\leq& \dpr(\overline{x_i},\zeta_i)
 +\dpr(\zeta_i,\zeta_{i+1}) +\dpr(\zeta_{i+1},\tilde{\zeta}_{i+1})\\
  &\leq& \frac{C\left(1+ \frac{1-\e}{2(1+\e)}\right)}
  {D^{3/2}\mun(q_i,\zeta_i)}
  +\frac{C\frac{1-\e}{12}}{D^{3/2}\mun(q_{i+1},\zeta_{i+1})}
  \quad \mbox{(by ({\bf a},$i$),  ({\bf d},$i$) and \eqref{eq:y})}\\
   &\leq& \frac{C(1+\e + \frac{7}{12} (1-\e))}
   {D^{3/2}\mun(q_{i+1},\zeta_{i+1})},
  \quad \mbox{(by (\ref{eq:b1}))}
\end{eqnarray*}
which proves ({\bf e},$i$).

Note that ({\bf x},$i$+1) and \eqref{eq:y}, together with
Proposition~\ref{prop:apps}, imply that $\mun(\tilde{q}_{i+1},
\tilde{\zeta}_{i+1})\leq  (1+\e)\mun(q_{i+1},\zeta_{i+1})$.
Also, that we have
$C(1+\e)(1+\e + \frac{7}{12} (1-\e))\le \nu_0\approx 0.3542$
and hence
$\dpr(\overline{x_i},\tilde{\zeta}_{i+1})\leq
\frac{\nu_0}{D^{3/2}
\mun(\tilde{q}_{i+1},\tilde{\zeta}_{i+1})}$.
We can therefore use Theorem~\ref{thm:gamma} to deduce
that $\overline{x_i}$ is an approximate zero
of~$\tilde{q}_{i+1}$ associated with its zero $\tilde{\z}_{i+1}$.
Therefore,
$x_{i+1}=N_{\tilde{q}_{i+1}}(\overline{x_i})$ satisfies
\begin{equation}%\label{eq:007}
  \dpr(x_{i+1},\tilde{\zeta}_{i+1})\leq
  \frac12\,\dpr(\overline{x_i},\tilde{\z}_{i+1})
  \leq \frac{C(1+\e + \frac{7}{12}(1-\e))}
  {2D^{3/2}\mun(q_{i+1},\zeta_{i+1})},
\end{equation}
where the last inequality is due to ({\bf e},$i$), and thus
\begin{eqnarray}\label{eq:007}
  \dpr(x_{i+1},\zeta_{i+1}) &\leq&
  \dpr(x_{i+1},\tilde{\zeta}_{i+1})
  + \dpr(\tilde{\zeta}_{i+1},\zeta_{i+1})\notag\\
  &\leq& \frac{\frac12 C(1+\e + \frac{7}{12}(1-\e)) +
  \frac{1}{12}C (1-\e)}{D^{3/2}\mun(q_{i+1},\zeta_{i+1})}\notag\\
  &=& \frac{C(\frac12 (1+\e) + \frac{3}{8}(1-\e))}
   {D^{3/2}\mun(q_{i+1},\zeta_{i+1})}.
\end{eqnarray}

Now we are ready to prove the last two implications.
We first show $({\bf{a}},i+1)$.

%From ({\bf c},$i$), ({\bf x},$i$) and (\ref{eq:b}) we have
%$$
% \dsp(q_i, \tilde{q}_{i+1})\leq \dsp(q_i,q_{i+1})+ \dsp(q_{i+1},
% \tilde{q}_{i+1})\leq \frac{\frac12 C(1-\e)
% +\frac{C(1-\e)}{12 (1+\e)}}{D^{3/2}\mun(q_i,\overline{x_i})},
%$$
%which together with Proposition~\ref{prop:apps} yields
%\begin{equation}\label{eq:011}
%  \frac{1}{1+\e}\mun(q_i,\overline{x_i})
%  \leq \mun(\tilde{q}_{i+1},\overline{x_i})
%  \leq (1+\e)\mun(q_i,\overline{x_i}).
%\end{equation}
%
Inequality \eqref{eq:bound_tildes} allows us to use once more
Proposition~\ref{prop:apps} to deduce
\begin{equation}\label{eq:011}
  \frac{1}{1+\e}\mun(\tilde q_i,\overline{x_i})
  \leq \mun(\tilde{q}_{i+1},\overline{x_i})
  \leq (1+\e)\mun(\tilde q_i,\overline{x_i}).
\end{equation}

Since $\overline{u}_i$ is less than or equal to
$\bfB(\tilde q_i,\overline{x_i})$ (by ({\bf u},$i$)),
from the choice of the
constant $k_2$ in Definition~\ref{def:k2}(ii) one has
\begin{eqnarray*}
    \overline{u}_i & \leq & \frac{\bfe}{(1+\e)^2 n D^{2}
    (\log N + D + n^2) \mun^2(\tilde q_i,\overline{x_i}) }\\
    & \leq & \frac{\bfe}{n D^{2} (\log N + D + n^2)
      \mun^2(\tilde q_{i+1},\overline{x_i}) }
      \quad \mbox{(by \eqref{eq:011})}.
\end{eqnarray*}
The condition on $u$ (for the pair
$(\tilde{q}_{i+1},\overline{x_{i}})$)
of Proposition~\ref{coro:Newton} is
thus verified, and applying this proposition we obtain
\begin{equation}\label{eq:008}
  \| \overline{x_{i+1}}-x_{i+1} \|
 = \Error (N_{\tilde{q}_{i+1}}(\overline{x_i})) \leq
  \frac{C (1-\e)}{4 \pi (1+\e)^2 D^{3/2}
  \mun(\tilde{q}_{i+1},\overline{x_i})}.
\end{equation}

%From ({\bf a},$i$) and Proposition~\ref{prop:apps} we have
%\begin{equation}\label{eq:005}
%\mun(q_i,\zeta_i)\leq (1+\e)\mun(q_{i},\overline{x}_i).
%\end{equation}
%therefore, using ({\bf a},$i$) \marginpar{\fbox{$C(1+\e)\leq C$}}
%$$
%\dsp(\overline{x_i},\zeta_i)\leq \frac{C}{D^{3/2}\mun(q_i, \zeta_i)}\leq \frac{(1+\e)C}{D^{3/2}\mun(q_{i+1}, \zeta_i)},
%$$
%which yields (again, by Proposition~\ref{prop:apps})
%\begin{equation}\label{eq:005}
%(1+\e)\mun(q_{i+1},\overline{x_i})\geq \mun(q_{i+1},\zeta_i).
%\end{equation}

The proof of \eqref{eq:bound_tildes} implicitly shows that
$\dsp(q_i,\tilde q_{i+1})\leq \frac{C(1-\e)}{D^{3/2}
\mun(q_{i},\z_{i})}$. Together with
({\bf a},$i$)
we are in the hypothesis of Proposition~\ref{prop:apps}
and we can deduce
$$
  \frac{1}{1+\e}\mun(q_i,\z_i)
  \leq \mun(\tilde{q}_{i+1},\overline{x_i})
  \leq (1+\e)\mun(q_i,\z_i).
$$
This inequality, together with
(\ref{eq:b1}), yields
\begin{equation}\label{eq:b2}
  \frac{1}{(1+\e)^2}\mun(q_{i+1},\z_{i+1})
  \leq \mun(\tilde{q}_{i+1},\overline{x_i})
  \leq (1+\e)^2\mun(q_{i+1},\z_{i+1})
\end{equation}
and using these bounds
%by combining \eqref{eq:b2} and \eqref{eq:011},
(\ref{eq:008}) becomes
\begin{equation}\label{eq:009}
   \| \overline{x_{i+1}}-x_{i+1} \| \leq
  \frac{C (1-\e)}{4\pi D^{3/2}\mun(q_{i+1},\zeta_{i+1})}.
\end{equation}
We now use this bound and the triangle inequality
to bound $\dpr(\overline{x_{i+1}},\zeta_{i+1})$ as follows
\begin{eqnarray*}
  \dpr(\overline{x_{i+1}},\zeta_{i+1})
    &\leq& \dpr(\overline{x_{i+1}},x_{i+1})
       +\dpr(x_{i+1},\zeta_{i+1}) \\
  &\leq& \frac{\pi}{2} \| \overline{x_{i+1}}-x_{i+1} \|
    +\dpr(x_{i+1},\zeta_{i+1})\\ %\quad\mbox{\tt caution here}\\
   &\leq& \frac{C(1-\e)}{8 D^{3/2}\mun(q_{i+1},\zeta_{i+1})}
    +\frac{C(1+\e + 3/4 (1-\e))}
     {2 D^{3/2}\mun(q_{i+1},\zeta_{i+1})}
   \quad \mbox{(by (\ref{eq:009}) and (\ref{eq:007}))}\\
   &=& \frac{C\left(\frac{1}{2}(1+\e) + \frac{3}{8} (1-\e)
     + \frac{1}{8} (1-\e)\right)}{D^{3/2}
     \mun(q_{i+1},\zeta_{i+1})} = \frac{C}{D^{3/2}
     \mun(q_{i+1},\zeta_{i+1})},
\end{eqnarray*}
which proves ({\bf a}) for $i+1$.

It remains to show $({\bf{u}},i+1)$. To do so note that
we may use ({\bf a},$i+1$) and ({\bf x}, $i+1$) together
with Proposition~\ref{prop:apps} to obtain
\eqref{eq:b} for $i+1$ (just as we obtained it for $i$).
Consequently,
\begin{eqnarray*}
\mun (\tilde{q}_{i+1}, \overline{x_{i+1}})
&\leq& (1+\e) \mun (q_{i+1}, \z_{i+1}) \\
&\leq& (1+\e)^2 \mun (q_{i}, \z_{i})
\qquad\mbox{(by (\ref{eq:b1}))}\\
&\leq& (1+\e)^3 \mun (\tilde q_{i}, \bar{x_{i}})
\qquad\mbox{(by (\ref{eq:b})).}
\end{eqnarray*}
Using this bound along with ({\bf u},$i$)  we obtain
\begin{eqnarray*}
     \overline{u}_i &\leq& \bfB(\tilde q_i,\overline{x_i})
     = \frac{k_2}{n D^{2} (\log N + D + n^2)
     \mun^2(\tilde q_i,\overline{x_i})}\\
    &\leq& \frac{k_2(1+\e)^6}{n D^{2} (\log N + D + n^2)
   \mun^2(\tilde q_{i+1},\overline{x_{i+1}})}
    = (1+\e)^6\bfB(\tilde{q}_{i+1},\overline{x_{i+1}}).
\end{eqnarray*}

We can therefore apply Proposition~\ref{prop:repeat_u}
with the pair $(\tilde{q}_{i+1},\overline{x_{i+1}})$
to deduce that
$\Error(\frac34\bfB(\tilde{q}_{i+1},\overline{x_{i+1}}))
\leq \frac14 \bfB(\tilde{q}_{i+1},\overline{x_{i+1}})$, and
consequently
\begin{eqnarray*}
   \left|\overline{u}_{i+1} - \frac{3}{4} \bfB(\tilde q_{i+1},
    \overline{x_{i+1}}) \right|
  &\leq& \frac{1}{4} \bfB(\tilde q_{i+1},\overline{x_{i+1}}),
\end{eqnarray*}
which proves ({\bf u},$i+1$).
\eproof

\subsection{Proof of Theorem~\ref{thm:main1}}

{\bf (i)\ } Since $(f,g,x,u)$ is an admissible input for
\ALHF\ we can use Proposition~\ref{prop:Recurence}
(and the notation therein). The estimate
$ \dpr(\overline{x_{k}},\zeta_{k}) \le
\frac{C}{D^{3/2}\mun(q_k,\zeta_k)}$ shown as
({\bf a},$k$) in that proposition
implies by Theorem~\ref{thm:gamma} that the
returned point~$\overline{x_k}$ is an approximate zero of
$q_k=f$ with associated zero~$\z_1$.
\medskip

{\bf (ii)\ }
The first instruction in \ALHF\ swaps $f$ by $-f$ if
$\dsp(\tilde f,g)\geq\frac{\pi}{2}$. The reason to do so
is that for nearly antipodal instances of $f$ and $g$
the difference $\dsp(f,\tilde f)$
may be arbitrarily magnified in $\dsp(q_\tau,\tilde q_\tau)$.
This does not occur under the assumption of infinite precision
and this is why such swap is not in the algorithms described
in~\cite{BePa08b,BC09}.

Let $h$ be either $-f$ or $f$ (according to whether \ALHF\
did the swap or not), $K(h,g,x)$ be the
number of iterations performed by \ALHF, and
$\{(q_\tau,\z_\tau)\}$ be the lifting of the path $E_{g,h}$.

Let $k\leq K(h,g,x)$ be a positive integer and consider any
$i\in\{0,\ldots,k-1\}$. Using
Proposition~\ref{prop:cdx} for $q_i,q_{i+1}$ together
with \eqref{eq:b} implies that, for all
$\tau\in[\bar{\tau_i},\bar{\tau_{i+1}}]$,
\begin{equation}\label{eq:b3}
   \frac{\mun(\tilde{q_i},\overline{x_i})}{(1+\e)^2}\leq
   \mun(q_\tau,\zeta_\tau)\leq (1+\e)^2
   \mun(\tilde{q_i},\overline{x_i}).
\end{equation}
Therefore,
\begin{eqnarray*}
  \int_{\overline{\tau_i}}^{\overline{\tau_{i+1}}}
   \mun^2(q_\tau,\zeta_\tau)d\tau
  &\geq&
  \int_{\overline{\tau_i}}^{\overline{\tau_{i+1}}}
  \frac{\mun^2(\tilde q_i,\overline{x_i})}{(1+\e)^4}d\tau\\
  &=&\frac{\mun^2(\tilde q_i,\overline{x_i})}{(1+\e)^4}
  (\overline{\tau_{i+1}}-\overline{\tau_{i}})\\
  &\geq& \frac{\mun^2(\tilde q_i,\overline{x_i})}{(1+\e)^4}
       \frac{3\lambda}{4\a D^{3/2}\mun^2(\tilde q_i,\overline{x_i})}
       \quad \mbox{(by Proposition~\ref{coro:DeltaTau})}\\
  &=&  \frac{3\lambda}{4(1+\e)^4\a D^{3/2}}.
\end{eqnarray*}
If $k=K(h,g,x)<\infty$ this implies
$$
  \int_0^1 \mun^2(q_\tau,\zeta_\tau)d\tau
   \ge \left( \frac{3\lambda}{4(1+\e)^4}\right)
   k\frac1{\a D^{3/2}} \geq k\frac1{408 \, \a D^{3/2}},
$$
which proves that
\begin{equation}\label{eq:boundK}
  K(h,g,x)\leq 408 \dsp(h,g) D^{3/2}
    \int_0^1 \mun^2(q_\tau,\zeta_\tau)d\tau = B(h,g,\z).
\end{equation}
It follows that the number of iterations
$K(f,g,x)$ satisfies either
$K(f,g,x)\leq B(f,g,\z)$ or
$K(f,g,x)\leq B(-f,g,\z)$.
Certainly ---and this introduces a factor of 2 but simplifies
the exposition---
\begin{equation*}
  K(f,g,x)\leq B(f,g,\z)+B(-f,g,\z).
\end{equation*}
In case $K(h,g,x)=\infty$ (a non-halting computation)
it implies that $\int_0^1 \mun^2(q_\tau,\zeta_\tau)d\tau=\infty$.

The bound for $\cost_{\ALHF}$ follows from the
$\Oh(N)$ cost of each iteration of \ALHF\ mentioned
in~ \S\ref{se:Newton}.
\medskip

{\bf (iii)\ }
For $i=1,\ldots,k-1$, due to ({\bf u},i),
$$
  \bar{u}_i\geq \frac{\bfB(\tilde q_i,\bar{x_i})}{2}
 =\Omega \left( \frac{1}{n D^{2} (\log N + D + n^2)
   \max_{ \tau \in [\bar{\tau_i},\bar{\tau_{i+1}}]}
   \mun^2(q_\tau,\zeta_\tau)} \right)
$$
the last by \eqref{eq:b3}. The statement now follows from
the equalities
\begin{equation}\tag*{\qed}
 u_*(f,g,\zeta)=\min_{i<k} \bar{u}_i
 \qquad\mbox{and}\qquad
 \mu_*(f,g,\zeta)=\max_{i<k}
    \max_{ \tau \in [\bar{\tau_i},\bar{\tau_{i+1}}]}
    \mun(q_\tau,\zeta_\tau).
\end{equation}

\section{Proof of Theorem~B}\label{sec:THB}

We follow here the proof of the corresponding
result for \MD\ in~\cite{BC09} and begin by recalling
two facts from this article. The first estimates
the mean square condition number on the path when
an extremity is fixed.

\begin{theorem}[Theorem 10.1 in~\cite{BC09}]\label{thm:IBA}
For $g\in S(\Hd)\setminus\Sigma$ we have
\begin{equation}\tag*{\qed}
   \E_{f\in S(\Hd)} \bigg(\dsp(f,g) \int_0^1 \mu_2^2(q_\tau) d\tau\bigg)
   \leq 818\,D^{3/2} N(n+1)\mum^2(g)+ 0.01 .
\end{equation}
\end{theorem}

The second bounds the condition of $\oU$.

\begin{lemma}[Lemma 10.5 in~\cite{BC09}]
\label{lemma:comp-munorm}
The maximum of the condition numbers $\mum(\oU) : =  \max_{\bz : \oU(\bz)=0} \{ \mun(\oU,\bz)\}$
satisfies
\begin{equation}\tag*{\qed}
\mum^2(\oU)   \le 2n\, \max_{i\leq n} \frac1{d_i} (n+1)^{d_i-1}\le 2\, (n+1)^{D}.
\end{equation}
\end{lemma}

The following proposition bounds the maximum
$\mu_*(f,g,\z)$ of the
condition number along a path from $(g,\z)$ to $f$
in terms of the number of iterations of \ALHF\ to follow this
path and of the condition number $\mun(g,\z)$ of the initial pair.

\begin{proposition}\label{prop:mu*}
Let $f,g\in S(\Hd)$ and $\z$ a zero of $g$.
The largest condition
number $\mu_*(f,g,\z)$ along the path
from $(g,\z)$ to $f$ satisfies
$$
  \mu_*(f,g,\z)\leq (1+\e)^{K(f,g,\z)}\mun(g,\z).
$$
\end{proposition}

\proof
Write $k:=K(f,g,\z)$ and let
$\mu_{*i}: = \displaystyle\max_{\tau\in [\overline{\tau_i},
\overline{\tau_{i+1}}]}\mun(q_\tau, \zeta_\tau)$. With this
notation, we have
$\mu_*(f,g,\z)=\displaystyle\max_{i=0,\ldots,k-1}\mu_{*i}$.
Furthermore, \eqref{eq:btau} states that, for all $i\leq k-1$,
$$
   \mu_{*i} \leq (1+\e) \mun(q_i, \z_i)
$$
and an immediate recursion yields
\begin{equation}\tag*{\qed}
 \mu_*(f,g,\z) = \max_{i\in \{1\ldots,k-1\}} \mu_{*i}
  \leq (1+\e)^{k}
 \mun(g,\z).
\end{equation}

We remark that
from the unitary invariance
of our setting,
for any unitary transformation $\nu \in \cU(n+1)$
and any $g\in\Hd$ and $x\in\proj^n$,
$$
\mun(g,x)=\mun(g\circ\nu^{-1},\nu x).
$$
Furthermore, for any execution of
\ALHF\ on an admissible input $(f,g,x)$, the number of iterations
$K(f,g,x)$ during the execution
satisfies $K(f,g,x) = K(f\circ\nu^{-1},g\circ\nu^{-1},\nu x)$
for any unitary transformation $\nu \in \cU(n+1)$.

But one can remark also that any zero $\bz_i$ of $\oU$ is the image
of $\bz_1= \frac{1}{\sqrt{2n}}(1,\ldots, 1)$
by a unitary transformation $\nu_i$ that leaves $\oU$ invariant.
%{\tt Prove this fact in a lemma?}
Thus, $K(f,\oU,\bz_1) = K(f\circ\nu_j^{-1},\oU, \bz_i)$
for all zeros $\bz_i$ of $\oU$, and
$\mum(\oU) = \mun(\oU,\bz_1)$.

We also obtain immediately

\begin{equation}\label{eq:KfUz1}
K(f,\oU,\bz_1) = \frac1{\Dn}\sum_{j=1}^{\Dn}
   K(f\circ\nu_j^{-1},\oU,\bz_j).
\end{equation}

But for all measurable functions $\varphi\colon S(\Hd)\to\R$
and all $\nu\in\cU(n+1)$ we have
\begin{equation*}\label{eq:inv1}
  \E_{f\in S(\Hd)}\varphi(f) =
  \E_{f\in S(\Hd)}\varphi(f\circ\nu),
\end{equation*}
due to the isotropy of the uniform measure on $S(\Hd)$.

Therefore,~\eqref{eq:KfUz1} implies
\begin{equation}\label{eq:ExpKfUz1}
\E_{f\in S(\Hd)} K(f,\oU,\bz_1) = \E_{f\in S(\Hd)} \frac1{\Dn}\sum_{j=1}^{\Dn}
   K(f,\oU,\bz_j).
\end{equation}

\proofof{Theorem~B}
From Lemma~\ref{lemma:comp-munorm},
$\mun(\oU,\bz_1) \leq \sqrt{2}\, (n+1)^{D/2}$, and thus
the initial value for $u$ in algorithm \MDF\ is
less than or equal to $\bfB(\oU,\bz_1)$.
Therefore, the tuple $(f,\oU,\bz_1,u)$ given to \ALHF\ during the
execution of \MD\ is an admissible input, and we can apply
Theorem~\ref{thm:main1} to that execution of \ALHF.
In particular, it follows that \MDF\ almost surely stops and
when it does so, it returns an approximate zero of $f$.
\smallskip

We next bound the average cost of \MDF.
Recall, we denoted by $K(f,\oU,\bz_1)$
the number of iterations of \ALHF\ during the
execution of \MDF\ with input $f$.
Again by Theorem~\ref{thm:main1},
for any root $\bz_j$ of $\oU$ we have
$$
  K(f,\oU,\bz_j) \le B(f, \oU, \bz_j) + B(-f, \oU, \bz_j).
$$
But we obviously have that $\E_{f\in S(\Hd)} B(f, \oU, \bz_j)= \E_{f\in S(\Hd)} B(-f, \oU, \bz_j)$,
and thus from~\eqref{eq:ExpKfUz1},
\begin{eqnarray*}
\E_{f\in S(\Hd)} K(f,\oU,\bz_1) & \leq &
2 \E_{f\in S(\Hd)} \frac{1}{\Dn} \sum_{j=1}^{\Dn} B(f, \oU, \bz_j) \\
     &  =   & 816 D^{3/2} \E_{f\in S(\Hd)}
        \dsp(f,\oU)     \int_0^1 \frac1{\Dn} \sum_{j=1}^{\Dn}
        \mun^2(q_\tau,\zeta_\tau^{(j)})d\tau \\
     & = & 816 D^{3/2} \E_{f\in S(\Hd)} \dsp(f,\oU)   \int_0^1
        \mu_2^2(q_\tau)d\tau,
\end{eqnarray*}
the last line by the definition of the mean square
condition number~\eqref{eq:mumean}.

Applying successively Theorem~\ref{thm:IBA} and
Lemma~\ref{lemma:comp-munorm}, we get
\begin{eqnarray}\label{eq:expK}
 \E_{f\in S(\Hd)} K(f,\oU,\bz_1)
& \leq &
   816 D^{3/2} (818\,D^{3/2} N(n+1)\mum^2(\oU)+ 0.01) \notag \\
& \leq & 816 D^{3/2} (818\,D^{3/2} N(n+1) \cdot 2\, (n+1)^{D} + 0.01) \notag \\
& = & 667488 D^3 N (n+1)^{D+1} + 8.16 D^{3/2} \notag \\
& \leq & 667489 \,D^3 N (n+1)^{D+1}.
\end{eqnarray}
The bound for the average of $\cost_{\MDF}$ follows from the
$\Oh(N)$ cost of each iteration of \ALHF.
\smallskip

We finally bound the average of the precision needed.
From Theorem~\ref{thm:main1} (iii), the finest precision
$u_*(f,\oU,\bz_1)$ along the execution of \ALHF\ (and therefore,
along that of \MDF\ ) satisfies, for some universal
constant $c$,
$$
  u_*(f,\oU,\bz_1)\geq
    \frac{1}{cn D^{2}
    (\log N + D + n^2)\mu_*^2(f,\oU,\bz_1)}.
$$
Hence by Proposition~\ref{prop:mu*},
\begin{eqnarray*}
  \big|\log u_*(f,\oU,\bz_1)\big| & \leq & \log\left(c n D^2 (\log N + D + n^2)
   \mu_*^2(f,\oU,\bz_1)\right)\\
 & = & 2\log\mu_*(f,\oU,\bz_1)+\log(c n D^2 (\log N + D + n^2))\\
  &\leq& 2K(f,\oU,\bz_1)\log (1+\e) +2\log \mun(\oU,\bz_1)
    +\log(c n D^2 (\log N + D + n^2)).
\end{eqnarray*}
%Using Proposition~\ref{prop:eq_dist}
%(now with $\varphi(\oU,\bz_1)=\log(1/u_*(f,\oU,\bz_1))$) we
Using Lemma~\ref{lemma:comp-munorm} and~\eqref{eq:expK},
we finally obtain
\begin{align*}
\E_{f\in S(\Hd)}
  \big|\log u_*(f,\oU,\bz_1)\big| \\
  \leq\; & 2 \E_{f\in S(\Hd)}
  \left(K(f,\oU,\bz_1)\log (1+\e) +\log \mun(\oU,\bz_1)
     \right.\\
   & +\,\left. \log(c n D^2 (\log N + D + n^2))\right)\\
    \leq\; & \log(1+\e) \cdot 1334978 D^3 N (n+1)^{D+1}\\
    & +\, D \log (\sqrt{2}(n+1)) +2 \log(c n D^2 (\log N + D + n^2))  \\
   =\; & \Oh(D^{3}N (n+1)^{D+1}).
\end{align*}
We observe that the initial precision
$u=\frac{k_2}{n D^2( \log N + D +n^2) 2 (n+1)^D}$
also satisfies this inequality.

\section{Proof of Theorem~A}\label{sec:THA}

We assume now a fixed round-off unit $\bu$.
In this context we consider the following trivial variation
of \ALHF:

\begin{center}
\algo
\> Algorithm \ALHFix\\[2pt]
\>{\bf input} $(f,g,x)$\\[2pt]
\>\> if $\dsp(f,g)\geq \frac{\pi}{2}$ then
 $g:=-g$\\[2pt]
\>\>$\tau:=0$, $q_\tau:=g$\\[2pt]
\>\> {\tt repeat}\\[2pt]
\>\>\>{\tt if $\bu>\frac34 \bfB(q_\tau,x)$ RETURN ``Failure''}\\[2pt]
\>\>\>$\Delta\tau:= \frac{\lambda}{\dsp(f,g) D^{3/2}\mun^2(q_\tau,x)}$\\[2pt]
\>\>\>$\tau:=\min\{1,\tau+\Delta\tau\}$\\[2pt]
\>\>\>$q_\tau:=t(\tau)f+(1-t(\tau)) g$\\[2pt]
\>\>\>$q_\tau:=\frac{q_\tau}{\|q_\tau\|}$\\[2pt]
\>\>\>$x:=N_{q_\tau}(x)$\\[2pt]
\>\>\>$x: = \frac{x}{\|x\|}$\\[2pt]
\>\> {\tt until} $\tau= 1$\\[2pt]
\>\> {\tt RETURN $x$}
\falgo
\end{center}
\medskip

The idea is simple: the precision in \ALHFix\ remains constant and
the algorithm proceeds until either it halts returning an approximate zero
$x$ of $f$ or it halts returning the message {\tt ``Failure''}. The latter
indicates that the precision is not sufficient to guarantee the correct
execution of the algorithm. In the former case, we say that \ALHFix\
{\em successfully halts}.

We can make \ALHFix\ and \ALHF\ even closer by taking advantage
of the level of generality we used to define the rounding functions
$r_u$ and the input-reading functions {\tt read\_input}. Recall,
the main property of $r_u$ is that $r_u(x)=x(1+\delta)$ with
$|\delta|\leq u$. Similarly, we
have made no assumptions on the functions
{\tt read\_input$_f$} besides the
fact that they return a rounded-off reading of the input system $f$
with the current precision.

For the rest of this section we assume that the rounding maps
$\{r_u\mid u\in (0,1)\}$ satisfy $r_u=r_{\bu}$ for all $u\geq \bu$.
We also assume that, for all $h\in\Hd$ the black-box
{\tt read\_input$_h$} is given by {\tt read\_input$_hh(\ )=r_u(h)$}
where $u$ is the current precision. All the results shown in
Sections~\ref{sec:homotopy} and \ref{sec:THB} hold in general and,
a fortiori, under these assumptions as well. In addition, we have
the following trivial lemma.

\begin{lemma}\label{lemma:key}
Let $h\in\Hd$ and $f=r_{\bu}(h)$.
If the computation of \ALHFix\ on input $(f,\oU,\bz_1)$ successfully
halts then this computation coincides with the computation of \ALHF\
on input $(h,\oU,\bz_1)$. In particular, they perform the same
number of iterations. Otherwise, both computations coincide until a
precision finer than $\bu$ is required, at which moment \ALHF\ proceeds
but \ALHFix\ halts with a failure message.\eproof
\end{lemma}
\medskip

\noindent
{\sc Proof of Theorem~A.\quad}
Algorithm \MDFix\ is what one would expect:

\algo
\> Algorithm \MDFix\\[2pt]
\>{\bf input} $f\in \Hu$ \\[2pt]
\>\> {\tt run \ALHFix\ on input $(f,\oU,\bz_1)$} \falgo
\medskip

Because of Lemma~\ref{lemma:key}, for $f\in\Hu$ we have
\begin{center}
\MDFix\ returns {\tt ``Failure''} with input $f$
$\iff$ for all $h\in r_{\bu}^{-1}(f)$, $u_*(h,\oU,\bz_1)<\bu$
\end{center}
where, we recall $u_*(h,\oU,\bz_1)$ is the finest $u_*$ required
by \MDF\  with input $h$.
Because of the equalities $\nu_{\bu}(\{f\})=\mu(r_{\bu}^{-1}(f))$
we therefore have
$$
  \Prob_{\nu_{\bu}}\{\mbox{\MDFix\ returns {\tt ``Failure''}}\}
  = \Prob_{h\sim N(0,\Id)}\{u_*(h,\oU,\bz_1)<\bu\}.
$$
Theorem~\ref{thm:main1}(iii) guarantees that, for some constant $A$,
$$
    u_*(h,\oU,\bz_1)\geq \frac{1}{AnD^2(\log N +D+n)\mu_*^2(f,\oU,\bz_1)}.
$$
Therefore,
\begin{align*}
\Prob_{\nu_{\bu}}&\{\mbox{\MDFix\ returns {\tt ``Failure''}}\}\\
%\Prob_{h\sim N(0,\Id)}&\{u_*(h,\oU,\bz_1)<\bu\} \\
\leq & \Prob_{h\sim N(0,\Id)}
\Big\{\mu_*^2(h,\oU,\bz_1) > \frac{1}{A\bu nD^2(\log N +D+n)}\Big\}\\
\leq& \Prob_{h\sim N(0,\Id)}
\Big\{ K(h,\oU,\bz_1) +\frac{\log\big(AnD^2(\log N +D+n)\mun^2(\oU,\bz_1)\big)}
{2\log(1+\varepsilon)}>
\frac{\log(1/\bu)}
{2\log(1+\varepsilon)}\Big\}
\end{align*}
the latter by Proposition~\ref{prop:mu*}. Let
$$
   X:=K(h,\oU,\bz_1) +\frac{\log\big(AnD^2(\log N +D+n)\mun^2(\oU,\bz_1)\big)}
{2\log(1+\varepsilon)}.
$$
Then $X$ is a positive random variable and
$$
  \E_{h\sim N(0,\Id)} X=\Oh(D^3N(n+1)^{D+1})
$$
by~\eqref{eq:expK} and the bound
$\mun^2(\oU,\bz_1)\leq 2(n+1)^D$ (Lemma~\ref{lemma:comp-munorm}).
We can now apply Markov's inequality to $X$ ---i.e.,
$\Prob\{X> t\}\leq \frac{\E X}{t}$--- and we finally
obtain
$$
\Prob_{\nu_{\bu}}\{\mbox{\MDFix\ returns {\tt ``Failure''}}\}
\,=\, \Oh\bigg(\frac{D^3N(n+1)^{D+1}}
 {\log(1/\bu)}\bigg).
$$
This shows the first assertion. To see the second,
let $\KFix(f,\oU,\bz_1)$ denote the number of iterations performed by
\ALHFix\ with input $(f,\oU,\bz_1)$. %Then
%Lemma~\ref{lemma:key} implies that,
%for all $f\in\Hu$ and all $h\in r_{\bu}^{-1}(f)$
%we have $\KFix(f,\oU,\bz_1)\leq K(h,\oU,\bz_1)$.
%The second assertion follows again from the equality
%$\nu_{\bu}(f)=\mu(r_{\bu}^{-1}(f))$ together
%with~\eqref{eq:expK}.

Let $\tau_H$ be the value of $\tau$ when $\MDFix$ halts on input
$(f,\oU,\bz_1),$ and let
$$
 \mu_{\bullet}(f,\oU,\bz_1) := \max_{\tau \in [0,\tau_H]}
 \left(\mun(q_\tau, \z_\tau)\right).
$$
By the proof of Theorem~\ref{thm:main1}(ii) and Lemma~\ref{lemma:key},
the number of iterations $\KFix(f,\oU,\bz_1)$ of $\MDFix$
is bounded as
\begin{equation}\label{eq:fix1}
\KFix(f,\oU,\bz_1) = \Oh\left( D^{3/2} \int_0^{\tau_H}
\mun^2 (q_\tau, \z_\tau) d\tau \right) = \Oh \left( D^{3/2}
\mu_{\bullet}^2(f,\oU,\bz_1)\right).
\end{equation}
Let $j$ be such that the maximum $\mu_{\bullet} (f,\oU,\bz_1)$ is attained in the interval $[\tau_j, \tau_{j+1}]$. From~\eqref{eq:btau},
\begin{equation}\label{eq:fix2}
 \mu_{\bullet}^2(f,\oU,\bz_1) \leq (1+\e)^2 \mun^2(q_j,\z_j).
% \leq \frac{(1+\e)^4 k_2}{n D^2 (\log N + D + n^2) \bu}.
\end{equation}
By Proposition~\ref{prop:Recurence}(u) and~\eqref{eq:b},
%at each step $i$ of the execution of $\MDFix,$
the precision $\bu$ satisfies
\begin{equation*}\label{eq:fix3}
\bu \leq \bfB(q,x)\leq (1+\e)^2 \bfB(q_j,\z_j)
= \frac{(1+\e)^2 k_2}{n D^2 (\log N + D + n^2) \mun^2(q_j, \z_j)}.
\end{equation*}
This inequality, together with~\eqref{eq:fix2}, imply
$$
  \mu_{\bullet}^2(f,\oU,\bz_1) \leq \frac{(1+\e)^4 k_2}
  {n D^2 (\log N + D + n^2) \bu}.
$$
and replacing this bound in~\eqref{eq:fix1} finally yields
\begin{equation}\tag*{\qed}
\KFix(f,\oU,\bz_1) = \Oh\left(\frac{1}{\sqrt{D} n (\log N + D + n^2) \bu}\right).
\end{equation}

{\small
%%  \bibliography{../../../book/book}

\begin{thebibliography}{10}

\bibitem{bast:83}
W.~Baur and V.~Strassen.
\newblock The complexity of partial derivatives.
\newblock {\em Theoretical Computer Science}, 22(3):317--330, 1983.

\bibitem{BeltranPardo08}
C.~Beltr\'an and L.~M. Pardo.
\newblock On {S}male's 17th problem: a probabilistic positive solution.
\newblock {\em Foundations of Computational Mathematics}, 8(1):1--43, 2008.

\bibitem{BePa08a}
C.~Beltr\'an and L.~M. Pardo.
\newblock Smale's 17th problem: average polynomial time to compute affine and
  projective solutions.
\newblock {\em Journal of American Mathematics Society}, 22(2):363--385, 2009.

\bibitem{BePa08b}
C.~Beltr\'an and L.~M. Pardo.
\newblock Fast linear homotopy to find approximate zeros of polynomial systems.
\newblock {\em Foundations of Computational Mathematics}, 11(1):95--129, 2011.

\bibitem{bcss:95}
L.~Blum, F.~Cucker, M.~Shub, and S.~Smale.
\newblock {\em Complexity and Real Computation}.
\newblock Springer-Verlag, 1998.

\bibitem{BC09}
P.~B\"urgisser and F.~Cucker.
\newblock On a problem posed by steve smale.
\newblock {\em Annals of Mathematics}.
\newblock To appear.

\bibitem{CKMM08}
F.~Cucker, T.~Krick, G.~Malajovich, and M.~Wschebor.
\newblock A numerical algorithm for zero counting, {I}: Complexity and
  accuracy.
\newblock {\em Journal of Complexity}, 24:582--605, 2008.

\bibitem{CS98}
F.~Cucker and S.~Smale.
\newblock Complexity estimates depending on condition and round-off error.
\newblock {\em Journal of the American Mathematics Society}, 46:113--184, 1999.

\bibitem{GoLoan}
G.~Golub and C.~Van~Loan.
\newblock {\em Matrix Computations}.
\newblock John Hopkins University Press, 3rd edition, 1996.

\bibitem{Higham96}
N.~Higham.
\newblock {\em Accuracy and Stability of Numerical Algorithms}.
\newblock SIAM, 96.

\bibitem{Shub93b}
M.~Shub.
\newblock Some remarks on {B\'e}zout's theorem and complexity theory.
\newblock In New~York Springer, editor, {\em From {T}opology to {C}omputation:
  {P}roceedings of the {S}malefest ({B}erkeley, {CA}, 1990)}, pages 443--455,
  1993.

\bibitem{Bez6}
M.~Shub.
\newblock Complexity of {B}\'ezout's theorem {VI}: {G}eodesics in the condition
  (number) metric.
\newblock {\em Foundations of Computational Mathematics}, 9(2):171--178, 2009.

\bibitem{Bez1}
M.~Shub and S.~Smale.
\newblock Complexity of {B}\'ezout's theorem {I}: {G}eometric aspects.
\newblock {\em Journal of the American Mathematics Society}, 6(2):459--501,
  1993.

\bibitem{Bez2}
M.~Shub and S.~Smale.
\newblock Complexity of {B}\'ezout's theorem {II}: {V}olumes and probabilities.
\newblock In F.~Eyssette and A.~Galligo, editors, {\em Computational Algebraic
  Geometry}, 109:265--285. Birkh\"auser, 1993.

\bibitem{Bez3}
M.~Shub and S.~Smale.
\newblock Complexity of {B}\'ezout's theorem {III}: {C}ondition number and
  packing.
\newblock {\em Journal of Complexity}, 9:4--14, 1993.

\bibitem{Bez5}
M.~Shub and S.~Smale.
\newblock Complexity of {B}\'ezout's theorem {V}: {P}olynomial time.
\newblock {\em Theoretical Computer Science}, 133:141--164, 1994.

\bibitem{Bez4}
M.~Shub and S.~Smale.
\newblock Complexity of {B}\'ezout's theorem {IV}: {P}robability of success;
  extensions.
\newblock {\em SIAM Journal of Numerical Analysis}, 33:128--148, 1996.

\bibitem{Smale86}
S.~Smale.
\newblock Newton's method estimates from data at one point.
\newblock In K.~Gross, R.~Ewing and C.~Martin, editors, {\em The Merging of
  Disciplines: New Directions in Pure, Applied, and Computational Mathematics},
  pages 265--285. Springer-Verlag, 1986.

\bibitem{smale:00}
S.~Smale.
\newblock Mathematical problems for the next century.
\newblock In {\em Mathematics: frontiers and perspectives}, pages 271--294.
  Amer. Math. Soc., Providence, RI, 2000.

\end{thebibliography}

}

\end{document}